\documentclass[12pt, reqno]{amsart}
\usepackage{latexsym}
\usepackage{amsmath}
\usepackage{amsthm}
\usepackage[mathscr]{eucal}
\usepackage{lmodern}
\usepackage[T1]{fontenc}
\usepackage{textcomp}
\usepackage[utf8]{inputenc}
\usepackage{amssymb}
\usepackage{enumerate}
\usepackage[abbrev]{amsrefs}
\usepackage{indentfirst}
\usepackage{graphicx}
\usepackage{bmpsize}
\usepackage{hyperref}
\usepackage{enumitem}
\usepackage[justification=justified]{caption}
\renewcommand{\thefootnote}{} 

\theoremstyle{plain} 
\newtheorem{theorem}{\indent\sc Theorem}[section]
\newtheorem{lemma}[theorem]{\indent\sc Lemma}
\newtheorem{corollary}[theorem]{\indent\sc Corollary}
\newtheorem{proposition}[theorem]{\indent\sc Proposition}

\newtheorem{subpath}{\indent\sc Subpath}
\theoremstyle{definition} 
\newtheorem{definition}[theorem]{\indent\sc Definition}
\newtheorem{remark}[theorem]{\indent\sc Remark}

\newenvironment{sproof}{%
  \proof}{\endproof}
\makeatletter
  
  \@addtoreset{equation}{section}
\makeatother
\allowdisplaybreaks 
\begin{document}


\title
[divergence function of the braided Thompson group]
{DIVERGENCE FUNCTION OF THE BRAIDED THOMPSON GROUP}
\author{YUYA KODAMA}
\date{}
\renewcommand{\thefootnote}{\arabic{footnote}}  
\setcounter{footnote}{0} 
\begin{abstract}
We prove that the braided Thompson group $BV$ has a linear divergence function. 
By the work of Dru{\c{t}}u, Mozes, and Sapir, 
this implies none of asymptotic cones of $BV$ has a cut-point. 
\end{abstract} 
\maketitle
\section{Introduction}
R. Thompson groups $F$, $T$, and $V$ are defined by Richard Thompson in 1965. 
These groups have many interesting properties. 
For instance, 
$F$ is the first example of a torsion-free group of 
type $F_\infty$ but not of type $F$, 
by Brown and Geoghegan \cite[Theorems 5.3 and 7.2]{brown1984infinite}. 
$T$ and $V$ are also group of type $F_\infty$, 
by Brown \cite[Theorem 4.17]{brown1987finiteness}, and
known as the first examples of infinite simple group 
with finite presentation, by Thompson. 
These groups have been studied using not only algebra but also analysis and geometry. 

On the other hand, various ``Thompson-like'' groups have been considered 
to study the relationship with Thompson groups and 
their own interesting properties. 
In this paper, we focus on the generalization of $V$, 
braided Thompson group $BV$ (sometimes this group is written as $V_{\mathrm{br}}$). 
This  group is defined independently 
by Brin \cite{brin2007algebra} and Dehornoy \cite{dehornoy2006group}. 
It is known that $BV$ has similar properties to those of $V$. 
For instance, Brin \cite[Theorem 5.1]{brin2006algebra} showed $BV$ is finitely presented, 
where the generators and relations are similar to those of $V$, and 
Bux, Fluch, Marschler, Witzel, and Zaremsky \cite[Main Theorem]{bux2016braided} proved that 
this group is also of type $F_\infty$. 
Zaremsky \cite{zaremsky2018geometric} suggests the relationship between 
$BV$ and metric spaces being $\mathrm{CAT(0)}$ or hyperbolic. 

Golan and Sapir \cite[Theorem 1.1]{golan2019divergence} showed that 
Thompson groups $F$, $T$, and $V$ have linear divergence functions. 
Roughly speaking, the divergence function of a finitely generated group $G$ is 
a function given by the length of the path connecting two points at the same distance from 
the origin while avoiding a small ball with the center at the origin in the Cayley graph. 
Gersten \cite[section 2]{gersten1994quadratic} introduced 
divergence of connected geodesic metric spaces as collections of such functions. 
We focus on each function rather than a collection, 
since it corresponds to the topological characterization of 
the asymptotic cones of the group \cite[Lemma 3.17]{dructu2010divergence}. 
Since braid groups have linear divergence functions (Proposition \ref{Bn_divergence}), 
it is natural to expect that so does braided Thompson group $BV$. 
In fact, Golan and Sapir posed a question whether their proof can be 
extended to Thompson-like groups. 
In this paper, we give a partial answer to this question. 
\begin{theorem} \label{main_theorem_BV}
Braided Thompson group $BV$ has a linear divergence function. 
\end{theorem} 

This paper is organized as follows. 
In Section \ref{background}, 
we summarize definitions of Thompson groups, braid groups, 
braided Thompson group, 
and in Section \ref{divergence_def}, we define the divergence functions of finitely generated groups. 
In Section \ref{main_proof}, 
first we prepare some lemmas on the number of carets of elements in $BV$. 
Then we construct a path which satisfies the requirement for the definition 
of the divergence function. 
This path is connecting two points $g$ in $BV$ and 
the point $v(|g|)$ in $F < BV$ which only depend the word length of $g$. 
This is achieved in the following way: 
For $g$, 
we construct the element $h$ (denoted by $w_1w_2w_3$ in Section \ref{main_proof}) 
in $BV$ such that $gh$ and $v(|g|)$ are commute. 
Then, we move $g \to gh \to ghv(|g|)=v(|g|)gh \to v(|g|)$. 
We remark that the above paths do not work for elements having less than three carets. For those elements $g$, we consider $gx_1$, a multiplication by a generator $x_1$, instead of $g$ itself. 

It is also interesting to study divergence functions of other Thompson-like groups, 
similar to $BV$. 
For example, Brady, Burillo, Cleary, and Stein \cite{brady2008pure} defined $BF$ (sometimes denoted by $F_{\mathrm{br}}$), 
which is braided version of Thompson group $F$. 
Acora and Cumplido \cite{aroca2020new} defined 
a family of infinitely braided Thompson's groups, 
which contains $BV$ as a special case. 
Another example 
is the Higman-Thompson groups, for instance. 
\subsection*{Acknowledgements}
I appreciate the referee for his or her close reading and precious comments. 
I would like to thank my supervisor, 
Professor Tomohiro Fukaya for his guidance. 
\section{Background} \label{background}
\subsection{Finitely generated groups and binary words}
A group $G$ is said to be \textit{finitely generated} 
if there exists a subset $X$ such that every element of $G$ 
can be written as a product of finitely many elements in $X \cup X^{-1}$, 
where $X^{-1}:=\{ x^{-1} \mid x \in X\}$. 
We call such a product a \textit{word} in $X$. 
We use ``$\equiv$'' and ``$=$'' 
to express equalities as words in $X$ and as elements of $G$, respectively. 
Let $x \equiv x_1x_2\cdots x_n$ be a word in $X$. 
A word $x^\prime$ is said to be a \textit{prefix} of $x$, denoted by $x^\prime \leq x$, 
if 
$x^{\prime} \equiv \emptyset$ or $x^\prime \equiv x_1 \cdots x_k$ for some $1 \leq k \leq n$, 
where $\emptyset$ denotes the empty word. 
A word $x^\prime$ is said to be a \textit{strict prefix} of $x$ if 
$x^\prime$ is a prefix of $x$ and $x^\prime \not \equiv x$. 

Let $w$ be a finite character string that consists of $0$ and $1$. 
We call such a character string \textit{binary word} and 
we also use ``$\equiv$'' to express equality. 
By the similar way above, we define a \textit{prefix} and \textit{strict prefix} of $w$.  
For every two binary words 
$w_1 \equiv a_1a_2\cdots a_j$ and $w_2 \equiv b_1b_2 \cdots b_k$ 
where $a_i$ and $b_i \in \{0, 1\}$ for every $i$, 
$w_1w_2$ denotes the concatenation $a_1a_2 \cdots a_j b_1b_2 \cdots b_k$. 
\subsection{Thompson groups}
\begin{figure}[tbp]
	\centering
	\includegraphics[width=40mm]{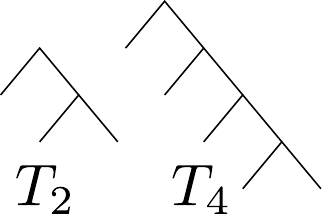}
	\caption{Examples of $T_n$}
	\label{binary_Tn}
\end{figure}

A \textit{rooted binary tree} is a tree with a distinguished vertex (\textit{root}) 
that has 2 edges, and vertices with either degree 1 (\textit{leaves}) or degree 3. 
We think of a rooted binary tree as a descending tree 
with the root as the only top vertex (level 0) and vertices of different levels. 
We define a \textit{caret} of a rooted binary tree to be a subtree of the tree 
that consists a vertex together with two downward-directed 
edge. 
We write \textit{all-right tree} $T_n$ for the rooted binary tree 
that is constructed by attaching a caret to the right edge of a caret $n$ times. 
Thus $T_1$ is a caret. 
See Figure \ref{binary_Tn}. 
The number of carets play important role to estimate the word lengths of elements of Thompson groups. 

\begin{figure}[tbp]
	\centering
	\includegraphics[width=130mm]{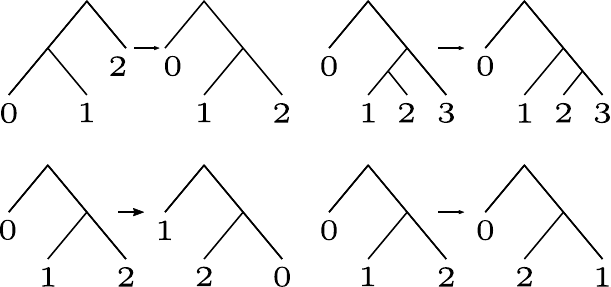}
	\caption{Examples of tree diagrams}
	\label{tree_diagrams}
\end{figure}

Let $(T_+, \sigma, T_-)$ be a triplet where 
$T_+$ and $T_-$ be finite rooted binary trees with $n$ caret, 
$L$ be the set of $(n+1)$ leaves and 
$\sigma$ be a permutation of $L$. 
We order the leaves of $T_+$ and $T_-$ from left to right from $0$ to $n$, 
respectively and 
use the numbers to represent the permutation $\sigma$. 
We call this \textit{tree diagram}. 
For example, see Figure \ref{tree_diagrams}. 

\begin{figure}[tbp] 
	\begin{tabular}{cc}
	\begin{minipage}{0.5\hsize}
		\centering
		\includegraphics[width=60mm]{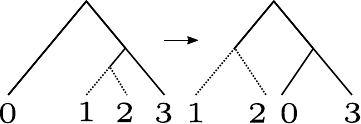}
		\captionsetup{justification=justified}
		\caption{Before the reduction}
		\label{not_reduced_V}
	\end{minipage}
	\begin{minipage}{0.5\hsize}
		\centering
		\includegraphics[width=60mm]{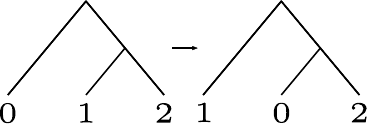}
		\caption{After the reduction and reordering the leaves}
		\label{reduced_V}
	\end{minipage}
	\end{tabular}
\end{figure}

Let $(T_+, \sigma, T_-)$ be the above tree diagram. We define a reduction of carets of a tree diagram as follows. We assume that two leaves $i$, $i+1$ have the same parent in $T_+$, two leaves $\sigma(i)$, $\sigma(i+1)$ have the same parent in $T_-$, and $\sigma(i+1)=\sigma(i)+1$ holds. In that case, each pair of the leaves forms carets. Then, we get the trees $T_+^\prime$ and $T_-^\prime$ by removing those carets. We regard the roots of the above carets as new leaves of the new trees, and we write $i_+$ and $i_-$ for the new leaves of $T_+^\prime$ and $T_-^\prime$, respectively. By sending $i_+$ to $i_-$ and sending other leaves by $\sigma$, we also get the permutation $\sigma^\prime$ on the set of $n$ leaves. This operation and its inverse are called \textit{reduction} and \textit{attachment} of carets, respectively. 
For example, see Figure \ref{not_reduced_V} and \ref{reduced_V}. 

Using these operations, we define the equivalence relation on the set of tree diagrams as follows. Two tree diagrams $(T_+, \sigma, T_-)$ and $(T_+^\prime, \sigma^\prime, T_-^\prime)$ are equivalent if $(T_+, \sigma, T_-)$ is obtained from $(T_+^\prime, \sigma^\prime, T_-^\prime)$ by a finite number of reductions and attachments. The Thompson group $V$ consists of all equivalence classes of tree diagrams. The product on $V$ is defined in the following way. 

\begin{figure}[tbp]
\centering
   \includegraphics[width=150mm]{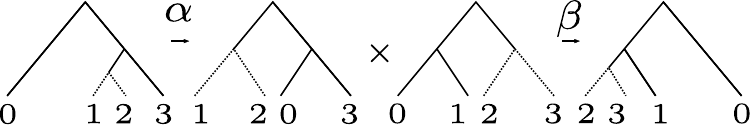}
    \caption{Diagrams of $a$ and $b$ with carets attached}
\label{g1g2keisan}
\end{figure}

\begin{figure}[tbp]
\centering
   \includegraphics[width=120mm]{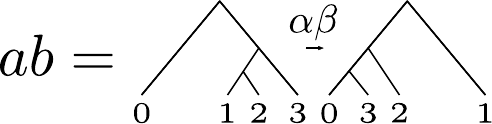}
    \caption{The product of the tree diagrams}
\label{g1g2kotae}
\end{figure}

For every two elements $a, b \in V$ represented by tree diagrams 
$(A_+, \alpha, A_-)$ and $(B_+, \beta, B_-)$, 
by successive attachments of carets, we get diagrams 
$(A_+^\prime, \alpha^\prime, A_-^\prime)$ and 
$(B_+^\prime, \beta^\prime, B_-^\prime)$ representing the same elements 
and such that $A_-^\prime=B_+^\prime$. 
Then the product $ab \in V$ is the equivalence class of 
$(A_+^\prime, \alpha^\prime\beta^\prime, B_-^\prime)$, 
where the permutation $\alpha^\prime\beta^\prime$ is composed from left to right. 
For example, see Figure \ref{g1g2keisan} and \ref{g1g2kotae}. 


The group $T$ is a subgroup of $V$ consists of equivalence classes of tree diagrams $(T_+, \sigma, T_-)$ where $\sigma$ is a cyclic permutation, and the group $F$ is a subgroup of $T$ consists of equivalence classes of tree diagrams $(T_+, \sigma, T_-)$ where $\sigma$ is the identity. 

For every caret of a rooted binary tree, we label its left edge by $0$ and the right edge by $1$. Since every leaf $o$ of such a tree $T$ corresponds to a unique path $s_T(o)$ 
from the root to the leaf, every leaf $o$ corresponds to a binary word $\ell_T(o)$ labeling the path from the root to $o$. We identify the path $s_T(o)$ with the binary word $\ell_T(o)$. 

By identifying the Cantor set $C$ with the set of infinite binary words, we can associate each tree diagram $(T_+, \sigma, T_-)$ to a homeomorphism from $C$ to itself. Indeed, for every leaf $o$ of $T_+$ and infinite binary word $w$, by mapping $\ell_{T_+}(o)w$ to $\ell_{T_-}(\sigma(o))w$, we get a homeomorphism. By the definition of this homomorphism $V \to \mathop{\mathrm{Homeo}}(C)$, the homeomorphism coming from $(T_+, \sigma, T_-)$ is the identity if and only if $\sigma$ is the identity and $T_+=T_-$,  so the homomorphism is injective. Hence, $V$ is a subgroup of the homeomorphism group of $C$. 

See \cite{cannon1996introductory} for details of the properties of Thompson groups. 
\subsection{Braid groups}
Let $n \in \mathbb{N}$. 
We briefly review the definition of geometric braid groups $B_n$. 
See \cite[Section 1.2]{kassel2008braid} for details. 
Let $I$ be the closed interval $[0, 1] \subset \mathbb{R}$. 
We call a topological space which is homeomorphic to $I$ \textit{topological interval}. 
\begin{definition}
[{\cite[Definition 1.4]{kassel2008braid}}]
A \textit{geometric braid} on $n$ strings is a set 
$b \subset \mathbb{R}^2 \times I$ formed 
by $n$ disjoint topological intervals called the \textit{strings} of $b$ such that 
the projection $\mathbb{R}^2 \times I \to I$ maps each string homeomorphically 
onto $I$ and 
\begin{align*}
b \cap (\mathbb{R}^2 \times \{ 0 \})= \{(0, 0, 0), (1, 0, 0), \dots, (n-1, 0, 0) \}, \\
b \cap (\mathbb{R}^2 \times \{ 1 \})= \{(0, 0, 1), (1, 0, 1), \dots, (n-1, 0, 1) \}. 
\end{align*}
We assume that every string goes from the bottom to up. 
\end{definition}
By the definition, every string of $b$ meets each plane 
$\mathbb{R}^2 \times \{t\}$ with $t \in I$ in exactly one point and 
connects a point $(i, 0, 0)$ to a point $(\sigma(i), 0, 1)$, where 
$\sigma$ is a permutation of $\{0, 1, \dots, n-1\}$. 
We call the both points \textit{endpoints} of the string, 
and call $\sigma$ \textit{underlying permutation} of the braid. 

\begin{definition}
[{\cite{kassel2008braid}}]
Two geometric braids $b$ and $b^\prime$ on $n$ strings are \textit{isotopic} 
if there exists a continuous map $F \colon b \times I \to \mathbb{R}^2 \times I$ such that 
for each $s \in I$, the map 
$F_s \colon b \to \mathbb{R}^2 \times I; x \mapsto F(x, s)$ is an embedding 
whose image is a geometric braid on $n$ strings, 
$F_0=\mathrm{Id} \colon b \to b$, and $F_1(b)=b^\prime$. 
Both the map $F$ and the family of geometric braids $\{F_s(b)\}_{s \in I}$ are called 
an \textit{isotopy} of $b$ to $b^\prime$. 
\end{definition}
The relation of isotopy is an equivalence relation 
on the class of geometric braids on $n$ strings. 
We call the equivalence classes and each string of an equivalence class 
\textit{braid} (on $n$ strands) and \textit{strand}, respectively. 
We write $B_n$ for the set of braids on $n$ strands. 

For every two geometric braids $b_1$ and $b_2$, 
we define their product $b_1b_2$ to be 
the set of points $(x, y, t) \in \mathbb{R}^2 \times I$ such that 
\begin{align*}
\mbox{$(x, y, 2t) \in b_1$ if $0 \leq t \leq \frac{1}{2}$},  
\end{align*}
and
\begin{align*}
\mbox{$(x, y, 2t-1) \in b_2$ if $\frac{1}{2} \leq t \leq 1$}. 
\end{align*}
It is clear that if $b_1$ and $b_2$ are isotopic to 
geometric braids $b_1^\prime$ and $b_2^\prime$, respectively, 
then the product $b_1b_2$ is isotopic to the product $b_1^\prime b_2^\prime$. 
Hence the product of $B_n$ is defined by the equivalence class of 
products of geometric braids. 

\begin{figure}[tbp]
	\centering
	\includegraphics[width=120mm]{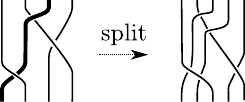}
	\caption{An example of splitting of the strand}
	\label{parallel_strand_Bn}
\end{figure}

A braid can be projected onto $\mathbb{R} \times \{0 \} \times I$ 
along the second coordinate with ``crossing information'' at each crossing point. 
Indeed, if necessary, by appropriate isotopies, 
we can assume that the number of strands involved in any intersection is two, 
every two strands meet transversely at each intersection point of the two strands, 
and there are only a finite number of such intersections. We call the intersections \textit{crossing points}. For each crossing point, the one with the lesser $y$-coordinate is denoted by \textit{over crossing}, and the other is denoted by the corresponding \textit{under crossing}. Then, we draw each over crossing by a continuous line, and each under crossing by a broken line. 
For example, Figure \ref{parallel_strand_Bn} are the projection of the elements in $B_4$ and $B_5$. 
In this paper, we identify braids with projected braids equipped with crossing information. 

We introduce an operation for braids which we use to define the product of elements of braided Thompson group. 

\begin{definition}[splitting] \label{splitting}
Let $0 \leq i \leq n-1$. 
Let $B$ be a braid on $n$ strands, 
$\{ b_k \mid k=0, 1, \dots, n-1 \}$ be the set of strands of $B$, 
$\sigma$ be the underlying permutation of $B$, and 
$(k, 0, 0)$ and $(\sigma(k), 0, 1)$ be endpoints of each strand $b_k$. 
We define a braid on $(n+1)$ strands $B^\prime$ to be the following: 
$B^\prime$ is obtained by adding strand $b_i^\prime$ to B 
such that it satisfies the following: 
\begin{enumerate}[font=\normalfont]
\item
Endpoints of $b_i^\prime$ are $(i+1/2, 0, 0)$ and $(\sigma(i)+1/2, 0, 0)$, 
then shift all endpoints appropriately so that they have integer $x$-values. 
\item
The strand $b_i^\prime$ does not cross with $b_i$. 
\item
The strand $b_i^\prime$ intersects with strands other than $b_i$ in the same way that $b_i$ intersects with strands in braid $B$. 
\end{enumerate}
In other words, $B^\prime$ is a braid 
such that $b_i^\prime$ is to the right of $b_i$ and 
the braid obtained from $B^\prime$ by removing $b_i$ is equal to $B$. 

We say that $b_i$ and $b_i^\prime$ are \textit{parallel}, and we call $B^\prime$ the \textit{splitting of the strand} $b_i$. 
\end{definition}
For example, see Figure \ref{parallel_strand_Bn}. 
\subsection{Braided Thompson group} 
Elements of Thompson groups $V$ can be 
seen as pairs of finite rooted binary trees, 
with permutations from leaves to itself. 
Roughly speaking, by replacing permutations with braids, we get elements of $BV$. 

\begin{figure}[tbp]
	\centering
	\includegraphics[width=60mm]{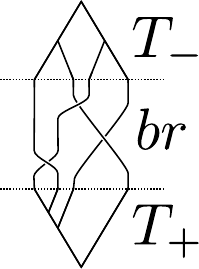}
	\caption{An example of tree-braid-tree diagram}
	\label{braided_example}
\end{figure}

Let $T_+$ and $T_-$ be finite rooted binary trees with $n$ carets and 
$br$ be a braid on $n+1$ strands from bottom to up. 
$(T_+, br, T_-)$ denotes a diagram where the leaves of both trees 
are joined by the braid with $T_+$ positioned upside down. 
We call this \textit{tree-braid-tree diagram}. 
For example, see Figure \ref{braided_example}. 

\begin{figure}[tbp]
	\centering
	\includegraphics[width=80mm]{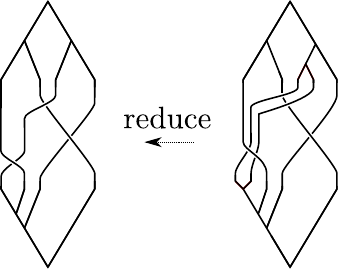}
	\caption{An example of reduction}
	\label{parallel_strand_BV}
\end{figure}

Let $(T_+, br, T_-)$ be the above tree-braid-tree diagram. Similar to Thompson groups, we define a reduction of carets of a tree-braid-tree diagram as follows. We assume that two strands $b_i$ and $b_i^\prime$ are parallel (cf.~Definition \ref{splitting}) and each endpoints of $b_i$ and $b_i^\prime$ have the same parent in $T_+$ and $T_-$. In that case, each pair of the endpoints (leaves) forms carets. Then, we get the trees $T_+^\prime$ and $T_-^\prime$ by removing those carets. We regard the roots of the above carets as new leaves of the new trees, and we write $i_+$ and $i_-$ for the new leaves of $T_+^\prime$ and $T_-^\prime$, respectively. By removing the strand $b_i^\prime$, letting the endpoints of $b_i$ be the new leaves, and keeping the other strands, we also get the braid $br^\prime$ on $n$ strands from $br$. This operation and its inverse operation are called \textit{reduction} of carets and \textit{splitting} of a strand, respectively. For example, see Figure \ref{parallel_strand_BV}. 

Using these operations, we define the equivalence relation on the set of tree-braid-tree diagrams as follows. Two tree-braid-tree diagrams $(T_+, br, T_-)$ and $(T_+^\prime, br^\prime, T_-^\prime)$ are equivalent if $(T_+, br, T_-)$ is obtained from $(T_+^\prime, br^\prime, T_-^\prime)$ by finite number of reductions and splittings. Each equivalence class has a unique representative with minimal number of carets. We call this diagram a \textit{reduced} tree-braid-tree diagram. 

The braided Thompson group $BV$ consists of all equivalence classes of tree-braid-tree diagrams. The product on $BV$ is defined in the following way. 

\begin{figure}[tbp]
	\centering
	\includegraphics[width=130mm]{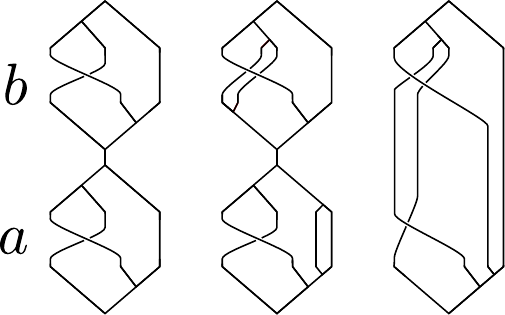}
	\caption{An example of the product $ab$ in $BV$}
	\label{operation_example}
\end{figure}

For every two elements $a, b \in BV$ represented by tree-braid-tree diagrams 
$(A_+, br_A, A_-)$ and $(B_+, br_B, B_-)$, by successive splittings of strands, 
we get diagrams $(A_+^\prime, br_A^\prime, A_-^\prime)$ and 
$(B_+^\prime, br_B^\prime, B_-^\prime)$ representing the same elements 
and such that $A_-^\prime = B_+^\prime$. 
Hence $br_A^\prime$ and $br_B^\prime$ are braids from the same braid group. 
Then the product $ab \in BV$ is 
the equivalence class of $(A_+^\prime, br_A^\prime br_B^\prime, B_-^\prime)$, 
where $br_A^\prime br_B^\prime$ is the braid that $br_A^\prime$ and $br_B^\prime$ 
connected from the bottom to the top, in this order. 
Figure \ref{operation_example} shows 
an example of a multiplication of elements of $BV$. 

\begin{figure}[tbp]
	\centering
	\includegraphics[width=130mm]{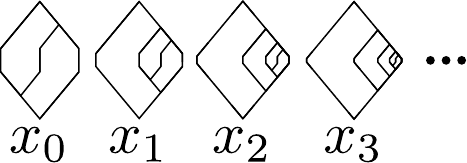}
	\caption{The infinite generators $x_i$}
	\label{generator_of_F_braided}
\end{figure}

\begin{figure}[tbp]
	\centering
	\includegraphics[width=135mm]{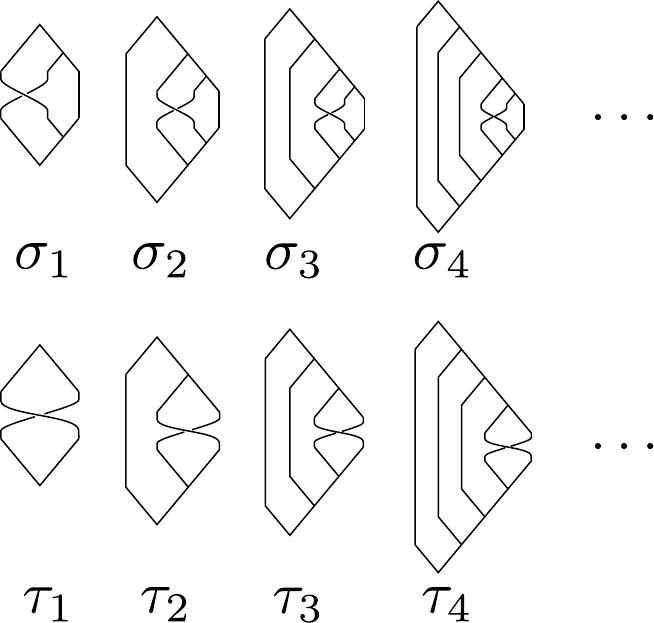}
	\caption{The infinite generators $\sigma_i$ and $\tau_i$}
	\label{generator_of_BV}
\end{figure}

It is known that $BV$ has the following infinite presentation. 
\begin{theorem}
[{\cite[Theorem 2.4]{brady2008pure}}] \label{presentation_of_BV}
The group $BV$ admits a presentation with generators: 
\begin{itemize}
\setlength{\parskip}{0cm}
\setlength{\itemsep}{0cm}
\item $x_i$, for $i \geq 0$, 
\item $\sigma_i$, for $i \geq 1$, 
\item $\tau_i$, for $i \geq 1$. 
\end{itemize}
and relators
\begin{description}
\setlength{\parskip}{0cm}
\setlength{\itemsep}{0cm}
\item[A] $x_jx_i=x_ix_{j+1}$, for $j>i$
\item[B1] $\sigma_i\sigma_j=\sigma_j\sigma_i$, for $j-i \geq 2$
\item[B2] $\sigma_i\sigma_{i+1}\sigma_i=\sigma_{i+1}\sigma_i\sigma_{i+1}$
\item[B3] $\sigma_i\tau_j=\tau_j\sigma_i$ for $j-i \geq 2$
\item[B4] $\sigma_i\tau_{i+1}\sigma_i=\tau_{i+1}\sigma_i\tau_{i+1}$
\item[C1] $\sigma_ix_j=x_j\sigma_i$, for $i < j$
\item[C2] $\sigma_ix_i=x_{i-1}\sigma_{i+1}\sigma_i$
\item[C3] $\sigma_ix_j=x_j\sigma_{i+1}$, for $i \geq j+2$
\item[C4] $\sigma_{i+1}x_i=x_{i+1}\sigma_{i+1}\sigma_{i+2}$
\item[D1] $\tau_ix_j=x_j\tau_{i+1}$, for $i-j \geq 2$
\item[D2] $\tau_ix_{i-1}=\sigma_i\tau_{i+1}$
\item[D3] $\tau_i=x_{i-1}\tau_{i+1}\sigma_i$. 
\end{description}
\end{theorem}
The reduced diagrams of generators $x_i$ are in Figure \ref{generator_of_F_braided}. 
The reduced diagrams of generators $\sigma_i$ and $\tau_i$ are in Figure \ref{generator_of_BV}. 
We note that a set of the generators $x_i$ corresponds to the standard infinite generating set of 
Thompson group $F$. 
Indeed, each $x_i$ is regarded as two rooted binary trees and 
identical permutation (see the upper low of Figure \ref{tree_diagrams}). 
Hence, $BV$ contains $F$ as a subgroup. 
Incidentally, in some papers, 
Thompson groups are defined by ``tree-permutation-tree diagrams'' 
similar to tree-braid-tree diagrams. 

Moreover, it is also known that $BV$ has the following finite presentation. 
\begin{theorem}
[{\cite[Theorem 3.1]{brady2008pure}}] \label{finite_presentation_of_BV}
The group $BV$ admits a finite presentation 
with generators $x_0, x_1, \sigma_1, \tau_1$ and relators
\begin{description}
\setlength{\parskip}{0cm}
\setlength{\itemsep}{0cm}
\item[a] $x_2x_0=x_0x_3$, $x_3x_1=x_1x_4$
\item[c1] $\sigma_1x_2=x_2\sigma_1$, $\sigma_1x_3=x_3\sigma_1$, $\sigma_2x_3=x_3\sigma_2$, $\sigma_2x_4=x_4\sigma_2$
\item[c3] $\sigma_2x_0=x_0\sigma_3$, $\sigma_3x_1=x_1\sigma_4$ 
\item[c4] $\sigma_1x_0=x_1\sigma_1\sigma_2$, $\sigma_2x_1=x_2\sigma_2\sigma_3$
\item[d1] $\tau_2x_0=x_0\tau_3$, $\tau_3x_1=x_1\tau_4$
\item[d2] $\tau_1x_0=\sigma_1\tau_2$, $\tau_2x_1=\sigma_2\tau_3$
\item[b1] $\sigma_1\sigma_3=\sigma_3\sigma_1$
\item[b2] $\sigma_1\sigma_2\sigma_1=\sigma_2\sigma_1\sigma_2$
\item[b3] $\sigma_1\tau_3=\tau_3\sigma_1$
\item[b4] $\sigma_1\tau_2\sigma_1=\tau_2\sigma_1\tau_2$
\end{description}
where the letters in the relators not in the set of $4$ generators are defined inductively 
by $x_{i+2}=x_i^{-1}x_{i+1}x_i$ for $i \geq 0$, $\sigma_{i+1}=x_{i-1}^{-1}\sigma_ix_i\sigma_i^{-1}$ for $i \geq 1$, and $\tau_{i+1}=x_{i-1}^{-1}\tau_i\sigma_i^{-1}$ for $i\geq 1$. 
\end{theorem}
We call $\{x_0, x_1, \sigma_1, \tau_1 \}$ \textit{standard generating set} of $BV$. 

As well as Thompson groups, for every caret of a rooted binary tree, 
we label its left edge by $0$ and the right edge by $1$. 
Since every leaf $o$ of such a tree $T$ corresponds to a unique path $s_T(o)$ 
from the root to the leaf, 
every leaf $o$ corresponds to a binary word $\ell_T(o)$ 
labeling the path from the root to $o$. 
We identify the path $s_T(o)$ with the word $\ell_T(o)$. 
The path $\ell_T(o)$ will be called a \textit{branch} of $T$. 
Let $(T_+, br, T_-)$ be a tree-braid-tree diagram of $g \in BV$, 
$o$ be a leaf of $T_+$, and $o^\prime$ be the corresponding leaf of $T_-$. 
We say that $\ell_{T_+}(o) \to \ell_{T_-}(o^\prime)$ is a \textit{branch} of 
the tree-braid-tree diagram $(T_+, br, T_-)$. 

Let $T$ be a rooted binary tree with $n$ carets. Recall that $T_n$ denotes an all-right tree. Then $(T, \mathrm{Id}, T_n) \in F$ is termed \textit{positive element}. Because there exist $0 \leq i_1 < i_2 < \dots < i_k$ and $r_1, r_2, \dots, r_k > 0$ such that
\begin{align*}
x_{i_1}^{r_1}x_{i_2}^{r_2} \cdots x_{i_k}^{r_k}=(T, \mathrm{Id}, T_n) 
\end{align*}
holds, where each $x_{i_t}$ is given by a diagram in Figure \ref{generator_of_F_braided} (See \cite[Theorem 2.5]{cannon1996introductory}). Similarly, $(T_n, \mathrm{Id}, T) \in F$ is termed \textit{negative element}. Since every element in $F$ is rewritten as a product of positive element and negative element, we call the product \textit{seminormal form}. For non-trivial element, let 
\begin{align*}
x_{i_1}^{r_1}x_{i_2}^{r_2} \cdots x_{i_k}^{r_k}x_{j_t}^{-s_t} \cdots x_{j_2}^{-s_2}x_{j_1}^{-s_1}
\end{align*} 
be a seminormal form, where $0 \leq i_1 < i_2 < \dots < i_k \neq j_t > \cdots j_2 >j_1 \geq 0$ and $r_1, r_2, \dots, r_k, s_1, \dots, s_t > 0$. This form is unique if we require the following condition: if $x_i$ and $x_i^{-1}$ exist in this form, then $x_{i+1}$ or $x_{i+1}^{-1}$ also exists. We call the unique form \textit{normal form}. By using the relation A, we can always get the normal form from the above seminormal form.  Furthermore, every such normal form represents non-identity element of $F$. See \cite[Corollary-Definition 2.7]{cannon1996introductory}. We call the part with positive exponents and  the one with negative exponents in the normal form \textit{positive part} and \textit{negative part}, respectively. 
\section{Divergence functions of finitely generated groups} \label{divergence_def}
Let $G$ be a finitely generated group, 
$X$ be a finite generating set of $G$, and 
$\Gamma$ be the Cayley graph $Cay(G, X)$. 
We will define the divergence functions of $G$. 
Since we consider asymptotic behavior of functions, we introduce 
a relation on the set of functions $\mathbb{R}_+ \to \mathbb{R}_+$ 
as follows. 
For such $f$ and $g$, we define $f \preceq g$ if 
\begin{align*}
f(x) \leq A g(Bx+C) + Dx + E
\end{align*}
for some $A, B, C, D, E \geq 0$ and all $x$. 
This defines an equivalence relation on the set of functions 
$\mathbb{R}_+ \to \mathbb{R}_+$, 
by saying $f \approx g$ if $f \preceq g$ and $g \preceq f$. 
We note that all linear functions and constant functions are equivalent. 

Let $\delta \in (0, 1)$. 
Then the \textit{$\delta$-divergence function} of $\Gamma$ is 
the smallest function $f_\delta(x)$ 
such that every two vertices of $\Gamma$ at distance $x$ from the identity can be 
connected by a path in $\Gamma$ of length less than 
$f_\delta (x)$ and 
avoiding the ball of radius $\delta x$ with a center at the identity. 
If no such path exists, take $f_\delta (x) = \infty$. 
For each $\delta \in (0, 1)$, 
the equivalence class of $f_\delta(x)$ is 
invariant under quasi-isometries, especially, it does not depend on the choice of finite generating set $X$. 
Hence the \textit{$\delta$-divergence function} of $G$ is 
defined to be the equivalence class of the $\delta$-divergence function of $\Gamma$. 

\begin{definition}
We say that the group $G$ has a \textit{linear divergence function} 
if there exists $\delta \in (0, 1)$ such that 
the $\delta$-divergence function of $G$
is equivalent to a linear function. 
\end{definition}

By definition, if $f_\delta$ is equivalent to a linear function, 
then for every $0 < \delta^\prime < \delta$, 
$f_{\delta^\prime}$ is equivalent to a linear function. 
Indeed, since a path that avoids the ball of radius $\delta x$ also avoids the ball of radius $\delta^\prime x$, 
$f_\delta(x) \geq f_{\delta^\prime}(x)$ holds for every $x$. 

Dru{\c{t}}u, Mozes and Sapir showed that having a linear divergence function is equivalent to the following topological property of asymptotic cones. 

\begin{theorem} [{\cite[correct version of Lemma 3.17]{dructu2010divergence,dructu2018corrigendum}}] \label{linear_cutpoint}
The following are equivalent. 
\begin{enumerate}[font=\normalfont]
\item 
$G$ has a linear divergence function. 
\item
For every $\delta \in (0, \frac{1}{54})$, 
$f_\delta$ is equivalent to a linear function. 
\item
None of asymptotic cones of $G$ has a cut-point. 
\end{enumerate}
\end{theorem}

We believe that the following result is well known. However, for reader's convenience, we give a sketch of a proof. 
\begin{proposition}\label{Bn_divergence}
For all $n \geq 3$, the braid group $B_n$ has a linear divergence function. 
\begin{sproof}
First, we note that the center of $B_n$ is isomorphic to $\mathbb{Z}$ (\cite[Theorem 1.24]{kassel2008braid}). 
Secondly, we also note that $B_n$ is not virtually cyclic, 
since $B_n$ has a subgroup which is isomorphic to $\mathbb{Z}^2$. 
Indeed, $B_3$ is a subgroup of $B_n$, and $B_3$ has a subgroup isomorphic to $\mathbb{Z}^2$ which is generated by commutative elements $p$ and $q$ (Figure \ref{p}, \ref{q}). 
\begin{figure}[tbp] 
	\begin{tabular}{cc}
	\begin{minipage}{0.5\hsize}
		\centering
		\includegraphics[height=50mm]{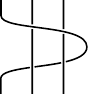}
		\captionsetup{justification=justified}
		\caption{the braid $p$}
		\label{p}
	\end{minipage}
	\begin{minipage}{0.5\hsize}
		\centering
		\includegraphics[height=50mm]{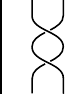}
		\caption{the braid $q$}
		\label{q}
	\end{minipage}
	\end{tabular}
\end{figure}
By combining the above two notes, we have that none of asymptotic cones of $B_n$ has a cut-point (\cite[Theorem 6.5]{dructu2005tree}). 
By Theorem \ref{linear_cutpoint}, this is equivalent to the property that $B_n$ has a linear divergence function. 
\end{sproof}
\end{proposition}
\begin{remark}
The above argument can work for pure braid groups as well. 
\end{remark}

\section{Proof of Theorem \ref{main_theorem_BV}} \label{main_proof}
\subsection{Number of carets for elements of $BV$}
Let $X=\{ x_0, x_1, \sigma_1, \tau_1 \}$ be the standard generating set of $BV$. 
For an element $g \in BV$, $|g|$ denotes the word length of $g$ with respect to the generating set $X$, 
and $N(g)$ denotes the number of carets in one of the trees 
in the reduced tree-braid-tree diagram of $g$. 
We will use the following estimate. 
\begin{theorem}[{\cite[Theorem 3.6]{burillo2009metric}}] \label{braided_bound}
For an element $g$ of $BV$ in tree-braid-tree diagram with $k$ total crossings, 
there exists a constant $C_1$ for which 
the word length satisfies the following inequalities: 
\begin{align*}
C_1 \max \{N(g), \sqrt[3]{k} \} \leq |g|. 
\end{align*}
\end{theorem}
Here we can assume that 
$0 < C_1 < 1$ without loss of generality. 

Let $g \in BV$ with a reduced tree-braid-tree diagram $(T_+(g), br(g), T_-(g))$. 
We call $T_+(g)$ the \textit{domain-tree} of $g$, 
$T_-(g)$ the \textit{range-tree} of $g$ and
$br(g)$ the \textit{braid} of $g$. 
Let $T$ be a rooted binary tree. 
Then, $\ell_0(T)$ denotes the length of the left most branch of T, 
that is, $\ell_0(T)=\ell$ if and only if $0^\ell$ is a branch of $T$, 
where we define 
\begin{align*}
i^\ell \equiv \underbrace{i \cdots i}_{\ell}, 
\end{align*}
for $i=0, 1$. 
Similarly, $\ell_1(T)$ denotes the length of the right most branch of T. 
For an element $g \in BV$, we define $\ell_i(g) := \ell_i(T_-(g))$, $i=0, 1$. 

We will need the following lemmas. 
Although the proofs of them 
are almost the same as in \cite{golan2019divergence}, 
we write down the proofs for reader's convenience. 
Note that the definition of $N(g)$ in this paper is 
different from that of $\mathcal{N}(g)$ in \cite{golan2019divergence}. 
The former denotes the number of carets 
and the latter the number of leaves, respectively. 

The following lemma corresponds to \cite[Lemma 2.2]{golan2019divergence}. 
\begin{lemma}\label{braided_lem_times_x0_2} 
Let $g$ be an element in $BV$ 
with reduced tree-braid-tree diagram $(T_+(g), br(g), T_-(g))$ 
and assume that $N(g) \geq 3$. 
Then
\begin{align*}
N(g) -1 \leq N(g x_0) \leq N(g)+1. 
\end{align*}
In addition, 
\begin{enumerate}[font=\normalfont]
\item
If $\ell_0(g)=1$ then $N(gx_0) = N(g) +1$ and 
$\ell_0(gx_0)=1$. 
\item
If $\ell_0(g) \neq 1$ then $N(gx_0) = N(g)$ or 
$N(gx_0) = N(g) -1$. 
Moreover, $\ell_0(gx_0) = \ell_0(g) - 1$. 
\item 
If $\ell_0(g) \neq 1$ and either $1$ or $01$ is 
a strict prefix of some branch of $T_-(g)$, 
then $N(gx_0) = N(g)$ and $1$ is a strict prefix of some branch of $T_-(gx_0)$, 
where $T_-(gx_0)$ is a range-tree of reduced tree-braid-tree diagram of $gx_0$. 
\end{enumerate}
\begin{proof}
We start by proving part (2). Assume that $\ell_0(g) \neq 1$. 
To multiply $g$ by $x_0$, we replace the reduced tree-braid-tree diagram 
of $x_0$ by an equivalent tree-braid-tree diagram $(R_+, \mathrm{Id}, R_-)$ where 
$R_+=T_-(g)$. 
Then $(T_+(g), br(g), R_-)$ is a tree-braid-tree diagram of the product $gx_0$. 
Let $u \to wv$ is a branch of $(T_+(g), br(g), T_-(g))$ 
where $w \equiv 00$, $w \equiv 01$, or $w \equiv 1$. 
By the construction of $(R_+, \mathrm{Id}, R_-)$, 
$u \to w^\prime v$ is a branch of $(T_+(g), br(g), R_-)$ 
where $w^\prime \equiv 0$, $w^\prime \equiv 10$, or $w^\prime \equiv 11$, respectively. 
All branches of $(T_+(g), br(g), R_-)$ can be written in this way. 

If $(T_+(g), br(g), R_-)$ is reduced diagram, then all assertions of part (2) hold, 
by the relation of $w \equiv 00$ and $w^\prime \equiv 0$. 
Indeed, it follows from the domain-tree that $N(gx_0)=N(g)$. Moreover, since the reduced diagram of $x_0$ has a branch $00 \to 0$, $\ell_0(R_-)=\ell_0(R_+)-1=\ell_0(T_-(g))-1=\ell_0(g)-1$. 
Hence, we can assume that $(T_+(g), br(g), R_-)$ is not reduced, that is, 
this diagram has a pair of branches $x0 \to y0$ and $x1 \to y1$ 
such that corresponding strands are parallel. 
Then $y \equiv 1$. 
Indeed, if $y$ is an empty word, $R_-$ has a branch $1$. 
This contradicts the fact that $\ell_1(x_0) \neq 1$. 
If $0$ is a prefix of $y$, then the diagram $(T_+(g), br(g), T_-(g))$ of $g$ 
has the pair of branches $x0 \to 0y0$ and $x1 \to 0y1$ 
such that corresponding strands are parallel, in contradiction to 
$(T_+(g), br(g), T_-(g))$ being reduced. 
If $10$ or $11$ is a prefix of $y$, 
the assumption that $(T_+(g), br(g), T_-(g))$ is reduced yields a 
contradiction in a similar way. 
Hence, $y \equiv 1$. 
So $R_-$ has branches of form $10$ and $11$. 
Now, we reduce the carets corresponding to $x0 \to 10$ and $x1 \to 11$ of 
the diagram $(T_+(g), br(g), R_-)$. 
Then the obtained diagram $(T_+^\prime(g), br^\prime(g), R_-^\prime)$ 
of $gx_0$ has a branch $x \to 1$ and this diagram is reduced. 
Indeed, if not reduced, there exists a pair of branches 
$x^\prime 0 \to y^\prime 0$ and $x^\prime 1 \to y^\prime 1$
such that corresponding strands are parallel. 
If $y^\prime$ is not an empty word, $0$ is a prefix of $y^\prime$. 
This contradicts with the same way as in branches $y0$ and $y1$. 
If $y^\prime$ is an empty word, then the obtained diagram has a branch 
$x^\prime 0 \to 0$. 
Since $N(g) \geq 3$ and $R_-$ has branches $10$ and $11$, 
$0$ is a strict prefix of some branch of $R_-$. 
Hence, $0$ is a strict prefix of some branch of $R_-^\prime$. 
This is a contradiction. 
Since the reduction is replacing branches $10$ and $11$ with $1$, 
we have $\ell_0(R_-^\prime)=\ell_0(R_-)$. 
Hence, part (2) holds. 

In the conditions of part (3) of the lemma, 
the tree-braid-tree diagram $(T_+(g), br(g), R_-)$ of $gx_0$ is reduced. 
Indeed, since $1$ or $01$ is a strict prefix of some branch of $T_-(g)$, 
by the relation of $w$ and $w^\prime$, 
either $11$ or $10$ is a strict prefix of some branch of $R_-$. 
It follows that $(T_+(g), br(g), R_-)$ is reduced, because, if not, as noted above, $R_-$ has branches both $11$ and $10$. 
In particular, $1$ is a strict prefix of some branch of $R_-$. 
Hence, part (3) holds. 

Now we assume the condition of part (1). 
To multiply $g$ by $x_0$, we replace the tree-braid-tree diagram $(T_+(g), br(g), T_-(g))$ 
by an equivalent diagram $(T_+^\prime(g), br^\prime(g), T_-^\prime(g))$ 
by attaching carets to the leaf of the branch $0$ of $T_-(g)$ and 
to the corresponding leaf of $T_+(g)$, 
and splitting of the corresponding strand. 
Let $(R_+, \mathrm{Id}, R_-)$ be the tree-braid-tree diagram of $x_0$ such that 
$R_+=T_-^\prime(g)$. 
Then we get the diagram $(T_+^\prime(g), br^\prime(g), R_-)$ of $gx_0$ and
we can proceed as part (2). 
To complete the proof it suffices to prove that 
$(T_+^\prime(g), br^\prime(g), R_-)$ is reduced. 
Indeed, in that case, by the construction, $N(gx_0)=N(g)+1$ and
since $\ell_0(R_+)=2$, we have $\ell_0(gx_0)=\ell_0(R_-)=1$. 
If $(T_+^\prime(g), br^\prime(g), R_-)$ is not reduced, 
there exists a pair of branches $x0 \to y0$ and $x1 \to y1$ such that 
corresponding strands are parallel. 
By the construction, the 
left most and second from the left 
branches of $T_-^\prime (g)$ are $00$ and $01$. 
Hence the left most and second from the left 
branches of $R_-$ are $0$ and $10$. 
This means $y$ is not empty and $0$ and $10$ are not prefixes of $y$. 
Then we first assume that $11$ is a prefix of $y$ and let $y \equiv 11y^\prime$ 
($y^\prime$ is probably an empty word). 
From the construction of $R_-$, 
$(T_+^\prime(g), br^\prime(g), T_-^\prime(g))$ has a pair of branch 
$x0 \to 1y^\prime0$ and $x1 \to 1y^\prime1$ such that 
corresponding strands are parallel. 
On the other  hand, $00$ and $01$ are only branches of $T_-^\prime(g)$ 
that can be reduced. 
This is a contradiction. 
Finally, we assume that $y \equiv 1$. Then, $(R_+, \mathrm{Id}, R_-)$ is the reduced tree-braid-tree diagram of $x_0$. Since $R_+=T_-^\prime(g)$ and 
$N(x_0)=2 $ (cf.~Figure \ref{generator_of_F_braided}) hold, by the construction of $T_-^\prime(g)$, $N(g)=2-1=1$. This contradicts the assumption of the lemma. 
Hence, part (1) holds. 
\end{proof}
\end{lemma}

The following corollary corresponds to \cite[Corollary 2.3]{golan2019divergence}. 
The proof given here is slightly modified. 
\begin{corollary} \label{braided_cor_times_x0_2}
Let $g$ be an element in $BV$ 
with a reduced tree-braid-tree diagram $(T_+(g), br(g), T_-(g))$ 
such that $N(g) \geq 3$. Let $\ell := \ell_0(g)$. 
Then the following assertions hold. 
\begin{enumerate}[font=\normalfont]
\item
If $N(g) \geq 3+(\ell-1)$ then for every $i \geq 0$ we have 
\begin{align*}
N(gx_0^i) \geq N(g) + i - 2(\ell - 1). 
\end{align*}
\item
If either $1$ or $01$ is a strict prefix of some branch of $T_-(g)$ 
then for every $i \geq 0$ we have
\begin{align*}
N(gx_0^i) = 
\max \{ N(g), N(g) + i - (\ell - 1) \}. 
\end{align*}
\end{enumerate}
\begin{proof}
To prove part (1), we first assume that $\ell=1$. 
By applying Lemma \ref{braided_lem_times_x0_2} (1) to $g$ 
iteratively, we have that for every $i \geq 0$, 
\begin{align*}
N(gx_0^i)=N(g)+i. 
\end{align*}
Thus, we can assume that $\ell >1$. 
Since $N(g) \geq 3 + (\ell-1)$, 
we can apply Lemma \ref{braided_lem_times_x0_2} (2) to $g$ at least $(\ell-1)$ times. 
Then we have $\ell_0(gx_0^{\ell-1})=1$ and for every $i \leq \ell-1$
\begin{align}
N(gx_0^i) \geq N(g)-i \geq N(g)-(\ell-1) \geq N(g) + i - 2(\ell-1). 
\label{braided_i_ell_1_1}
\end{align}
Since $N(gx_0^{\ell-1}) \geq 3$ and $\ell_0(gx_0^{\ell-1})=1$, 
by applying Lemma \ref{braided_lem_times_x0_2} (1) to $gx_0^{\ell-1}$ iteratively, 
we have 
\begin{align}
N(gx_0^{\ell-1}x_0^j)=N(gx_0^{\ell-1})+j 
\geq N(g)-(\ell-1)+j, \label{braided_Substituting}
\end{align}
 for every $j \geq 0$. 
By substituting $j=i-(\ell-1)$ in inequality \eqref{braided_Substituting}, 
we have that for every $i \geq \ell-1$, 
\begin{align}
N(gx_0^i) \geq N(g) + i -2(\ell-1). \label{braided_i_ell_1_2}
\end{align}
It follows from inequalities \eqref{braided_i_ell_1_1} and \eqref{braided_i_ell_1_2} that 
for every $i \geq 0$, 
\begin{align*}
N(gx_0^i) \geq N(g)+i-2(\ell-1), 
\end{align*}
as required. 

In the condition of part (2), we first assume that $\ell=1$, again. 
By Lemma \ref{braided_lem_times_x0_2} (1) (applying iteratively), 
we have
\begin{align*}
N(gx_0^i)=N(g)+i=\max \{N(g), N(g) + i\}, 
\end{align*}
for every $i \geq 0$. 
Thus, we can assume that $\ell >1$. 
Since $N(g) \geq 3$ and either $1$ or $01$ is a strict prefix of some branch of $T_-(g)$, 
by Lemma \ref{braided_lem_times_x0_2} (2) and (3), 
$N(gx_0^i)=N(g) \geq N(g)+i-(\ell-1)$ for every $i \in \{ 0, \dots, \ell-1 \}$ 
and $\ell_0(gx_0^{\ell-1})=1$. 
Thus, it suffices to prove that for every $i \geq \ell-1$ we have 
$N(gx_0^i)=N(g)+i-(\ell-1) \geq N(g)$. 
Since $\ell_0(gx_0^{\ell-1})=1$ and $N(gx_0^{\ell-1})=N(g) \geq 3$, 
by Lemma \ref{braided_lem_times_x0_2} (1), we have 
\begin{align}
N(gx_0^{\ell-1}x_0^j)=N(gx_0^{\ell-1})+j=N(g)+j, \label{braided_i_ell_1_3}
\end{align}
for every $j \geq 0$. 
Substituting $j=i-(\ell-1)$ in inequality \eqref{braided_i_ell_1_3} gives that 
for every $i \geq \ell-1$, 
\begin{align*}
N(gx_0^i)=N(g)+i-(\ell-1), 
\end{align*}
as required. 
\end{proof}
\end{corollary}

The proofs of the following lemma and corollary are symmetric to
those of Lemma \ref{braided_lem_times_x0_2} 
and Corollary \ref{braided_cor_times_x0_2}. 
We only need to switch $0$ and $1$. 
\begin{lemma}\label{braided_lem_times_x0^-1_2}
Let $g$ be an element in $BV$ with reduced tree-braid-tree diagram $(T_+(g), br(g), T_-(g))$ 
and assume that $N(g) \geq 3$. 
Then
\begin{align*}
N(g) -1 \leq N(g x_0^{-1} )\leq N(g)+1. 
\end{align*}
In addition, 
\begin{enumerate}[font=\normalfont]
\item 
If $\ell_1(g)=1$ then $N(gx_0^{-1}) = N(g) +1$ and 
$\ell_1(gx_0^{-1})=1$. 
\item 
If $\ell_1(g) \neq 1$ then $N(gx_0^{-1}) = N(g)$ or 
$N(gx_0^{-1}) = N(g) -1$. 
Moreover, $\ell_1(gx_0^{-1}) = \ell_1(g) - 1$. 
\item 
If $\ell_1(g) \neq 1$ and either $0$ or $10$ is 
a strict prefix of some branch of $T_-(g)$, 
then $N(gx_0^{-1}) = N(g)$ and $0$ is 
a strict prefix of some branch of $T_-(gx_0^{-1})$, 
where $T_-(gx_0^{-1})$ is a range-tree of 
reduced tree-braid-tree diagram of $gx_0^{-1}$. 
\end{enumerate}
\end{lemma}
\begin{corollary} \label{braided_cor_times_x0^-1_2}
Let $g$ be an element in $BV$ with reduced tree-braid-tree diagram $(T_+(g), br(g), T_-(g))$ 
such that $N(g) \geq 3$. Let $\ell := \ell_1(g)$. 
Then the following assertions hold. 
\begin{enumerate}[font=\normalfont]
\item 
If $N(g) \geq 3+(\ell-1)$ then for every $i \geq 0$ we have 
\begin{align*}
N(gx_0^{-i}) \geq N(g) + i - 2(\ell - 1). 
\end{align*}
\item
If either $0$ or $10$ is a strict prefix of some branch of $T_-(g)$ 
then for every $i \geq 0$ we have
\begin{align*}
N(gx_0^{-i}) = 
\max \{ N(g), N(g) + i - (\ell - 1) \}. 
\end{align*}
\end{enumerate}
\end{corollary}

\begin{figure}[tbp]
	\centering
	\includegraphics[width=45mm]{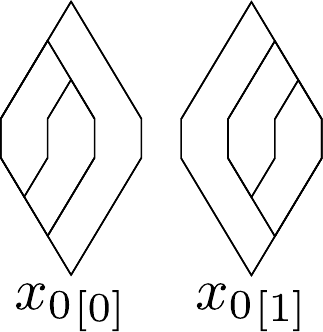}
	\caption{The reduced tree-braid-tree diagrams of ${x_0}_{[0]}$ and ${x_0}_{[1]}$}
	\label{braided_copy_x0}
\end{figure}

\begin{figure}[tbp]
	\centering
	\includegraphics[width=60mm]{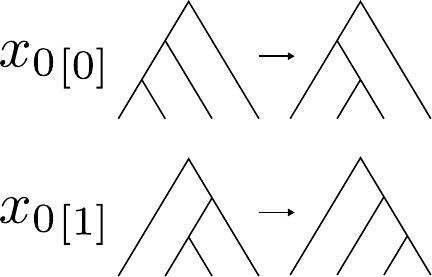}
	\caption{The reduced tree diagrams of ${x_0}_{[0]}$ and ${x_0}_{[1]}$}
	\label{braided_copy_x0_triple}
\end{figure}

\begin{figure}[tbp]
	\centering
	\includegraphics[width=140mm]{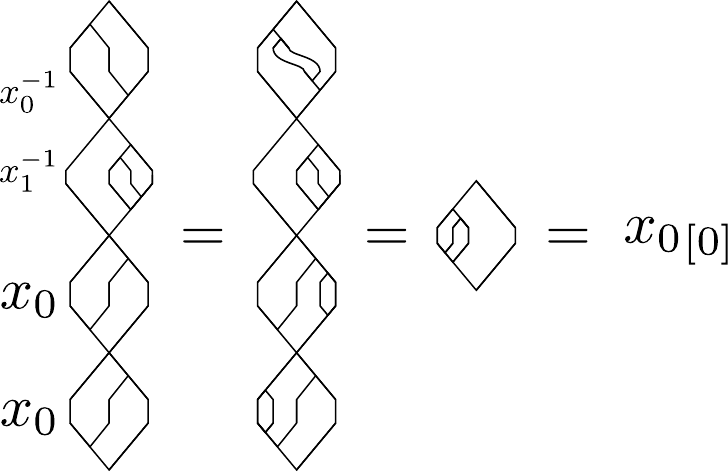}
	\caption{Calculation of $x_0^2x_1^{-1}x_0^{-1}={x_0}_{[0]}$}
	\label{x0^2x1^-1x0^-1=x0_0_braided}
\end{figure}

The next lemma describes the result of multiplying 
an element of $BV$ on the right by 
an element of $F$ with a following specific form. 
Let $u$ be a finite binary non-empty word 
and $h \in F$ be an non-identity element 
with reduced tree-braid-tree diagram $(T_+(h), \mathrm{Id}, T_-(h))$. 
Let $T$ be a minimal finite rooted binary tree which contains the branch $u$. 
We take two copies of the tree $T$. 
To the first copy, we attach the tree $T_+(h)$ at the end of the branch $u$, 
and we write $R_+$ for this tree. 
To the second copy, we attach the tree $T_-(h)$ at the end of the branch $u$, 
and we write $R_-$ for this tree. 
Then the element $h_{[u]}$ is the one represented by the tree-braid-tree diagram, 
where domain-tree is $R_+$, range-tree is $R_-$ 
and braid is the ``identity'', that is, all strands are straight. 
It is clear from the definition that $h_{[u]} \in F < BV$. 
For example, ${x_0}_{[0]}$ and ${x_0}_{[1]}$ are 
elements corresponding to the diagrams 
in Figure \ref{braided_copy_x0} or \ref{braided_copy_x0_triple}. 
Note that $x_0^2x_1^{-1}x_0^{-1}={x_0}_{[0]}$ holds, 
see Figure \ref{x0^2x1^-1x0^-1=x0_0_braided}. 

The following lemma corresponds to \cite[Lemma 2.6]{golan2019divergence}. 
\begin{lemma}\label{braided_add_copy}
Let $g \in BV$ be a non-identity element, 
$u \to v$ be a branch of $g$, 
$h$ be a non-identity element of $F$. 
Let $h^{\prime} = h_{[v]}$. 
Then 
\begin{align*}
N(gh^{\prime}) = N(g) + N(h). 
\end{align*}
Moreover, if h consists of branches $w_i \to z_i, i=1, \dots, k$ and 
$B$ is the set of branches of $g$ which are not equal to $u \to v$, 
then $gh^{\prime}$ consists of branches $uw_i \to vz_i, i=1, \dots, k$ 
along with all branches in $B$. 
\begin{proof}
Let $(T_+(g), br(g), T_-(g))$, $(T_+(h), \mathrm{Id}, T_-(h))$ 
and $(T_+(h^\prime), \mathrm{Id}, T_-(h^\prime))$ be the 
reduced tree-braid-tree diagrams of $g$, $h$ and $h^\prime$, respectively. 
To multiply $g$ by $h^\prime$, we note that 
the minimal refinement of $T_-(g)$ and $T_+(h^\prime)$ is the tree 
obtained from $T_-(g)$ by attaching the tree $T_+(h)$ at the bottom of 
the branch $v$, 
since $T_-(g)$ has a branch $v$ and $T_+(h^\prime)$ is constructed from the minimal tree which has a branch $v$. 
Let $S$ denote the described tree 
and $(R_1, br^\prime(g), S)$ be an equivalent tree-braid-tree diagram of $g$. 
We note that 
$R_1$ is obtained from $T_+(g)$ by attaching a copy of $T_+(h)$ 
to the bottom of the branch $u$. 
If $(S, \mathrm{Id}, R_2)$ is a tree-braid-tree diagram of $h^\prime$, 
we also note that 
$R_2$ can be obtained from $T_-(g)$ by attaching $T_-(h)$ 
at the bottom of the branch $v$. 
The product of the tree-braid-tree diagrams 
$(R_1, br^\prime(g), S)$ and $(S, \mathrm{Id}, R_2)$ is $(R_1, br^\prime(g), R_2)$. 
Since $(R_1, br^\prime(g), S)$ has branches $uw_i \to vw_i$ and branches in $B$
($x \to y$ denotes these one), and 
$(S, \mathrm{Id}, R_2)$ has branches $vw_i \to vz_i$ and $y \to y$, 
it follows that $(R_1, br^\prime(g), R_2)$ has branches $uw_i \to vz_i$ and $x \to y$. 
To finish the proof, it remains to prove that 
$(R_1, br^\prime(g), R_2)$ is reduced. 
Since $(T_+(h), \mathrm{Id}, T_-(h))$ is reduced, 
$(T_+(h), \mathrm{Id}, T_-(h))$ has no pair of branches 
of the form $p0 \to q0$ and $p1 \to q1$ 
where $p0, p1 \equiv w_i$ for some $i$, respectively and 
$q0, q1 \equiv z_i$ for some $i$, respectively. 
Hence, $(R_1, br^\prime(g), R_2)$ has no pair of branches of the form 
$up0 \to vq0$ and $up1 \to vq1$, that is, 
$(R_1, br^\prime(g), R_2)$ has no pair of the branches of the form $uz_i \to vq_i$ 
such that reducible. 
Similarly, since $(T_+(g), br(g), T_-(g))$ is reduced, 
$(R_1, br^\prime(g), R_2)$ has no pair of branches of the form $x \to y$ such that reducible. 

Recall that $R_1$ is obtained from $T_+(g)$ by attaching a copy of $T_+(h)$. 
Then it is clear that $N(gh^\prime)=N(g)+N(h)$ holds, and the proof is complete. 
\end{proof}
\end{lemma}

\subsection{Construction of the path}
If $w$ is a word over the alphabet $X$, $\| w \|$ denotes the length of $w$. 
Note that any word $w$ over the alphabet $X$ can be regarded 
as an element of $BV$, 
then we have $|w| \leq \| w \|$. 

\begin{remark}
Golan-Sapir constructed a path between elements whose number of carets is greater than or equal to three in \cite[Proposition 2.7]{golan2019divergence}. Linear divergence of Thompson groups $F$, $T$, $V$ follow immediately from this path. However, in the case of the braided Thompson group $BV$, we need a little more discussion. Because the number of $g \in BV$ such that $N(g) \leq 2$ is infinite. For example, $\tau_1, \tau_1^2, \tau_1^3, \dots$ all have one caret. 
\end{remark}
First, we consider elements in $BV$ whose number of carets are greater than or equal to three (Proposition \ref{braided_main_prop}). Next, we construct paths between elements whose number of carets are less than three and others. 

The following proposition corresponds to 
\cite[Proposition 2.7]{golan2019divergence}. 
In \cite{golan2019divergence}, they constructed the path that consists of five subpaths, subpath 1, \dots, subpath 5. In this paper, we will take a similar process, but 
our subpath 3 (and therefore also the path $w$) is different from the original one. 
Our subpath 3 does not work for Thompson group $T$, 
but an almost similar approach works for Thompson groups $F$ and $V$. 
\begin{proposition} \label{braided_main_prop}
There exist constants $\delta, D > 0$ 
and a positive integer $Q$ such that the following holds. 
Let $g \in BV$ be an element with $N(g) \geq 3$. 
Then there exists a path of length at most $D|g|$ in the Cayley graph 
$\Gamma = \mathrm{Cay}(BV, X)$ which avoids a $\delta|g|$-neighborhood 
of the identity and which has initial vertex $g$ 
and terminal vertex $x_0^{Q|g|}x_1^{-1}x_0^{-Q|g|+1}$. 

In other words, there exists a word $w$ in the alphabet $X$ 
such that $ \|w \| < D|g|$; for any prefix $w^{\prime}$ of $w$, 
we have $|gw^{\prime}| > \delta|g|$ and such that 
\begin{align*}
gw = x_0^{Q|g|}x_1^{-1}x_0^{-Q|g|+1}. 
\end{align*}
\begin{proof}
Let $C_1$ be the constant from Theorem \ref{braided_bound}. 
We give 5 subwords $w_1, \dots, w_5$ and then let $w \equiv w_1 \cdots w_5$. 
Let $(T_+(g), br(g), T_-(g))$ be the reduced tree-braid-tree diagram of $g$. 
\begin{subpath}
If $0$ is not a branch of $T_-(g)$ we let $w_1 \equiv \emptyset $ and let $g_1 = g$. 
Otherwise, we let $w_1 \equiv x_0^2x_1^{-1}x_0^{-1}$ and let $g_1 = gw_1$. 
\end{subpath}
Let $(T_+(g_1), br_1, T_-(g_1))$ be the reduced tree-braid-tree diagram of $g_1$. 
The following lemma corresponds to \cite[Lemma 2.8]{golan2019divergence}. 
\begin{lemma}\label{braided_subpath1}
We have that $0$ is not a branch of $T_-(g_1)$. 
Moreover, $N(g) \leq N(g_1) \leq N(g)+2$ hold, and for every prefix 
$w^{\prime}$ of $w_1$, we have $N(gw^{\prime}) \geq N(g)$. 
\begin{proof}
If $0$ is not a branch of $T_-(g)$ then $g_1=g$ and $w_1\equiv\emptyset$, so the lemma holds. 
Thus, we can assume that $0$ is a branch of $T_-(g)$. 
Let $u$ be the binary word such that $(T_+(g), br(g), T_-(g))$ has the branch $u \to 0$. 
We recall that $w_1=x_0^2x_1^{-1}x_0^{-1}={x_0}_{[0]}$ 
(cf.~Figure  \ref{x0^2x1^-1x0^-1=x0_0_braided}). 
Hence, by Lemma \ref{braided_add_copy}, 
$uv_1 \to 0v_2$ is a branch of 
the reduced tree-braid-tree diagram of 
$g_1=g w_1=g{x_0}_{[0]}$ 
for each branch $v_1 \to v_2$ of $x_0$. 
Therefore, $0$ is not a branch of $T_-(g_1)$ 
since it is a strict prefix of some branch. 

For the second claim, 
by Lemma \ref{braided_add_copy}, 
we have 
\begin{align*}
N(g_1)=N(g{x_0}_{[0]})=N(g)+N(x_0)=N(g)+2. 
\end{align*}

For the last claim, we will consider the number of carets of 
$gx_0$, $gx_0^2$ and $gx_0^2x_1^{-1}$. 
Since $0$ is a branch of $T_-(g)$, we have $\ell_0(g)=\ell_0(T_-(g))=1$. 
Hence, by Lemma \ref{braided_lem_times_x0_2} (1), 
$N(gx_0)=N(g)+1$ and $\ell_0(gx_0)=1$. 
Again, by applying Lemma \ref{braided_lem_times_x0_2} (1) to $gx_0$, 
$N(gx_0^2)=N(gx_0)+1=N(g)+2$. 
Finally, we note that $gx_0^2x_1^{-1}=g_1x_0$ and $N(g_1)=N(g)+2$. 
By applying the inequality in Lemma \ref{braided_lem_times_x0_2} to $g_1$, 
we have 
\begin{align*}
N(gx_0^2x_1^{-1})=N(g_1x_0) \geq N(g_1)-1=N(g)+1, 
\end{align*}
and the proof is complete. 
\end{proof}
\end{lemma}
\begin{subpath}
We fix an integer $M \geq 100/C_1$. 
Then we define a word $w_2$ by
\begin{align*}
w_2 \equiv x_0^{-M(N(g_1)+1)}x_1x_0^{M(N(g_1)+1)}
\end{align*}
and we let $g_2 = g_1w_2$. 
\end{subpath}
Let $(T_+(g_2), br_2, T_-(g_2))$ be the reduced tree-braid-tree diagram of $g_2$. 
The following lemma corresponds to \cite[Lemma 2.9]{golan2019divergence}. 
\begin{lemma}\label{braided_subpath2}
The following assertions hold. 
\begin{enumerate}[font=\normalfont]
\item For every prefix $w^{\prime}$ of $w_2$, we have 
$N(g_1w^{\prime}) \geq N(g_1)$. 
\item $N(g_2) \geq M N(g_1)$. 
\end{enumerate}
\begin{proof}
We first prove part (2). 
It follows from the relation A in Theorem \ref{presentation_of_BV} that, 
as an element of $BV$, we have $w_2=x_m$ where $m=M(N(g_1)+1)+1$. 
Let $\ell_1=\ell_1(g_1)$ and $u$ be a finite binary word 
such that $u \to 1^{\ell_1}$ is a branch of $(T_+(g_1), br_1, T_-(g_1))$. 
By considering the minimum tree-braid-tree diagram 
where some branch is $1^{\ell_1}$, 
$\ell_1 \leq N(g_1)$ is clear. 
We note that $x_m={x_{m-\ell_1}}_{[1^{\ell_1}]}$. 
Hence, by Lemma \ref{braided_add_copy}, we have
\begin{align*}
N(g_2)=N(g_1w_2)=N(g_1x_m)=N(g_1{x_{m-\ell_1}}_{[1^{\ell_1}]})
=N(g_1)+N(x_{m-\ell_1}). 
\end{align*}
By the definition of standard infinite generating set of $F$, 
we have $N(x_{m-\ell_1})=m-\ell_1+2$. Hence, 
\begin{align*}
N(g_2)=N(g_1)+m-\ell_1+2
=N(g_1)-\ell_1+M(N(g_1)+1)+3
\geq MN(g_1), 
\end{align*}
as $\ell_1 \leq N(g_1)$. 
Thus, part (2) holds. 

Before proceeding the proof of part (1), we note that if $z_j \to q_j$, $j=1, \dots, n$ 
are the branches of $x_{m-\ell_1}$, then, by Lemma \ref{braided_add_copy}, 
the branches of the diagram $(T_+(g_2), br_2, T_-(g_2))$ of 
$g_2=g_1x_m$ are $uz_j \to 1^{\ell_1}q_j$, $j=1, \dots, n$ 
as well as all branches $a_k \to b_k$ 
which are branches of the diagram $(T_+(g_1), br_1, T_-(g_1))$, 
other than $u \to 1^{\ell_1}$. 
By Lemma \ref{braided_subpath1}, 
0 is a strict prefix of some branch of $T_-(g_1)$. 
Hence, $0$ is a strict prefix of some branch of $T_-(g_2)$. 

Now, let $w^{\prime}$ be a prefix of $w_2$. 
Then either (a) $w^\prime \equiv x_0^{-i}$ for some $0 \leq i \leq M(N(g_1)+1)$, 
or (b) $w^\prime \equiv x_0^{-M(N(g_1)+1)}x_1x_0^i$ 
for some $0 \leq i \leq M(N(g_1)+1)$. 

By Lemma \ref{braided_subpath1}, $N(g_1) \geq N(g)$ and 
$0$ is a strict prefix of some branch of $T_-(g_1)$. 
Hence, by Corollary \ref{braided_cor_times_x0^-1_2} (2), 
for $g_1$ and any $i \geq 0$, $N(g_1x_0^{-i}) \geq N(g_1)$. 
Hence, part (1) of the lemma holds for 
prefixes $w^\prime$ of type (a). 
Next, we consider the element $g_1w^\prime$ for $w^\prime$ of type (b). 
As an element of $BV$, we have
\begin{align*}
g_1w^\prime &=g_1x_0^{-M(N(g_1)+1)}x_1x_0^i \\
&=g_1x_0^{-M(N(g_1)+1)}x_1x_0^{M(N(g_1)+1)} \cdot x_0^{i-M(N(g_1)+1)} \\
&=g_1x_m x_0^{i-M(N(g_1)+1)} \\
&=g_2x_0^{i-M(N(g_1)+1)}. 
\end{align*}
From the note above, $0$ is a strict prefix of some branch of $T_-(g_2)$, 
and we have already shown that $N(g_2) \geq MN(g_1) \geq N(g)$. 
Hence, by Corollary \ref{braided_cor_times_x0^-1_2} (2), 
for $g_2$ and any $s \geq 0$, we have
\begin{align*}
N(g_2x_0^{-s}) \geq N(g_2) \geq MN(g_1). 
\end{align*}
We also note that $i-M(N(g_1)+1) \leq 0$. 
By substituting $-s=i-M(N(g_1)+1)$, 
we have 
\begin{align*}
N(g_1w^\prime)=N(g_2x_0^{i-M(N(g_1)+1)}) \geq MN(g_1), 
\end{align*}
as required. 
\end{proof}
\end{lemma}
\begin{subpath}
For reduced tree-braid-tree diagram $(T_+(g_1), br_1, T_-(g_1))$ of $g_1$, 
we assume that some strand of $br_1$ connects $0$th leaf to $k$th leaf, 
where $0 \leq k \leq N(g_1)$. 
Let $h$ be the element of $BV$ given by the (maybe not reduced) tree-braid-tree diagram 
$(T_-(g_1), br_h, T_+(g_1))$, where $br_h$ is following. 
If $k > 0$, the $k$th strand goes over 
the $(k-1)$-th strand, $(k-2)$-th strand, \dots, $0$th strand, in order, 
and other strands are straight. 
If $k=0$, all strands are straight. 
For example, Figure $\ref{w3_example}$ illustrates the construction of braid 
for $N(g_1)=5$ and $k=4$. 
Now, we let $w_3$ be the minimal word over $X$ such that $h=w_3$ in $BV$ 
and let $g_3=g_2w_3$. 
\end{subpath}
\begin{figure}[tbp]
	\centering
	\includegraphics[width=80mm]{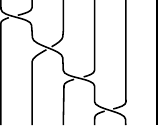}
	\caption{The braid of $h$ for $N(g_1)=5$ and $k=4$}
	\label{w3_example}
\end{figure}
Part (2) and (3) of the following lemma correspond to \cite[Lemma 2.10]{golan2019divergence}. 
\begin{lemma} \label{braided_subpath3}
The following assertions hold. 
\begin{enumerate}[font=\normalfont]
\item $\| w_3 \| \leq 14 N(g_1) $. 
\item $N(g_3) \geq (M-1)N(g_1)+M+3$. 
\item $\ell_0(g_3) \leq N(g_1)+1$ and 
$0^{\ell_0(g_3)} \to 0^{\ell_0(g_3)}$ is a branch of $g_3$. 
\end{enumerate}
\begin{proof}
Let $n=N(g_1)$. 
To prove part (1), we split $(T_-(g_1), br_h, T_+(g_1))$ into 
three diagrams by using all-right trees $T_n$, 
giving a split into a positive element $(T_-(g_1), \mathrm{Id}, T_n) \in F$, 
a braid element $(T_n, br_h, T_n)$, 
and a negative element $(T_n, \mathrm{Id}, T_+(g_1)) \in F$, 
where $T_n$ has $n$ carets (recall Figure \ref{binary_Tn}). 
Let $p$, $Br_h$ and $q$ be the minimal words over X such that 
$p=(T_-(g_1), \mathrm{Id}, T_n)$, 
$Br_h=(T_n, br_h, T_n)$ and 
$q=(T_n, \mathrm{Id}, T_+(g_1))$ in $BV$, respectively. 
We identify these words with elements of $BV$. 

First, we prove that $|p| \leq 6n$ holds. 
Let $\mathcal{A}$ be the standard generating set of $F$ 
such that $\mathcal{A} \subset X$ 
and we note that $|p|_\mathcal{A} \leq 6N(p)$
(see \cite[Theorem 1 and Proposition 2]{burillo2001metrics}). 
We also note that the tree-braid-tree diagram $(T_-(g_1), \mathrm{Id}, T_n)$ might not be reduced. 
Hence, we have 
\begin{align*}
|p| \leq |p|_{\mathcal{A}} \leq 6N(p) \leq 6n. 
\end{align*}
Similarly, we have $|q| \leq 6n$. 

Next, we prove that $|Br_h| \leq 2n$ holds. 
To get this upper bound, we rewrite the word $Br_h$ 
by elements of infinite generator of $BV$. 
The following rewritings are obvious. 
\begin{alignat*}{2}
k&=n & &\Rightarrow Br_h=\tau_n \sigma_{n-1} \dots \sigma_1, \\
k&=n-1 & &\Rightarrow Br_h=\sigma_{n-1} \dots \sigma_1, \\
k&=n-2 & &\Rightarrow Br_h=\sigma_{n-2} \dots \sigma_1, \\
&\; & &\;\; \vdots \\
k&=1 & &\Rightarrow Br_h= \sigma_1, \\
k&=0 & &\Rightarrow Br_h=\emptyset. 
\end{alignat*}
It suffices to consider only the case $k=n$, 
as we can get the following estimation. 
Indeed, for $n \geq 4$, we have
\begin{align*}
Br_h 
&= \tau_n \sigma_{n-1} \cdots \sigma_2 \sigma_1 \\
&= (x_0^{-(n-2)} \tau_2 x_0^{n-2})(x_0^{-(n-3)} \sigma_2 x_0^{n-3})
\cdots (x_0^{-1}\sigma_2 x_0) ( \sigma_2 \sigma_1) \\
&= (x_0^{-(n-2)}\tau_2x_0)(\sigma_2 x_0) 
\cdots (\sigma_2 x_0) (\sigma_2 \sigma_1) \\
&= (x_0^{-(n-2)} \sigma_1^{-1}\tau_1 x_0 x_0)
(x_0^{-1}\sigma_1x_1 \sigma_1^{-1} x_0)
\cdots (x_0^{-1} \sigma_1 x_1 \sigma_1^{-1} x_0)
(x_0^{-1} \sigma_1 x_1 \sigma_1^{-1} \sigma_1) \\
&= (x_0^{-(n-2)} \sigma_1^{-1} \tau_1 x_0)
(\sigma_1 x_1 \sigma_1^{-1})
\cdots (\sigma_1 x_1 \sigma_1^{-1})(\sigma_1 x_1) \\
&= (x_0^{-(n-2)} \sigma_1^{-1} \tau_1 x_0)\sigma_1 x_1^{n-2}, 
\end{align*}
where we rewrite $\tau_n=x_0^{-(n-2)}\tau_2x_0^{n-2}$, $\sigma_i=x_0^{-(i-2)}\sigma_2x_0^{i-2}$ for each $i \geq 3$, $\tau_2=\sigma_1^{-1}\tau_1x_0$, and $\sigma_2=x_0^{-1}\sigma_1x_1\sigma_1^{-1}$ by the relations D1($j=0$), C3($j=0$), D2($i=1$), and C2($i=1$) in Theorem \ref{presentation_of_BV}, respectively. 
Hence, we have
\begin{align*}
|Br_h| \leq n-2+3 + 1+ n-2 = 2n. 
\end{align*}
When $n=3$, we have 
\begin{align*}
Br_h &= \tau_3 \sigma_2 \sigma_1 \\
&= (x_0^{-1} \tau_2 x_0)(\sigma_2 \sigma_1) \\
&= (x_0^{-1} \sigma_1^{-1} \tau_1 x_0 x_0)(x_0^{-1} \sigma_1 x_1) \\
&= x_0^{-1} \sigma_1^{-1} \tau_1 x_0 \sigma_1 x_1. 
\end{align*}
Hence, we have
\begin{align*}
|Br_h| \leq 6 = 2 \times 3 =2n. 
\end{align*}
Therefore, we have
\begin{align*}
\| w_3 \| =|h| \leq |p| + |Br_h| + |q| \leq 14n, 
\end{align*}
as required. 

For part (2) and (3), we recall the 
form of branches of 
the reduced tree-braid-tree diagram $(T_+(g_2), br_2, T_-(g_2))$ of $g_2$ 
as described in the proof of Lemma \ref{braided_subpath2}. 
Let $\ell_1=\ell_1(T_-(g_1))$ and let $m=M(N(g_1)+1)+1$. 
Then 
\begin{align*}
g_2=g_1 x_m=g_1 {x_{m-\ell_1}}_{[1^{\ell_1}]}. 
\end{align*}
Now, let $u$ be such that 
$u \to 1^{\ell_1}$ is a branch of $(T_+(g_1), br_1, T_-(g_1))$. 
If $z _j \to q_ j$, $j = 1, \dots, n$ are the branches of 
reduced tree-braid-tree diagram of $x_{m-\ell_1}$ 
then the branches of 
$(T_+(g_2), br_2, T_-(g_2))$ are $u z_ j \to 1^{\ell_1}q_j$, $j = 1, \dots, n$ 
as well as all the branches $a_k \to b_k$ of $(T_+(g_1), br_1, T_-(g_1))$, 
other than $u \to 1^{\ell_1}$. 

Let $v$ be such that $1^{\ell_1} \to v$ is a branch of the tree-braid-tree diagram
$(T_-(g_1), br_h, T_+(g_1))$ of $h$. 
Then $uz_j \to vq_j$, $j=1, \dots, n$ are all branch of 
reduced tree-braid-tree diagram of $g_2h$. 
Indeed, by the proof of Lemma \ref{braided_add_copy}, 
$T_-(g_2)$ is constructed by attaching 
the range-tree $T_-(x_{m-\ell_1})$ of 
reduced tree-braid-tree diagram of $x_{m-\ell_1}$ 
at the end of branch $1^{\ell_1}$ of the $T_-(g_1)$, 
and so $T_-(g_1)$ is a rooted subtree of $T_-(g_2)$. 
Then to multiply $g_2$ by $h$, 
we replace $(T_-(g_1), br_h, T_+(g_1))$ 
by an equivalent tree-braid-tree diagram $(T_+(h), br_h^\prime, T_-(h))$
where $T_+(h)=T_-(g_2)$. 
By construction, $1^{\ell_1} q_j \to v q_j$ are 
branches of $(T_+(h), br_h^\prime, T_-(h))$. 
Hence, 
\begin{align*}
(T_+(g_2), br_2, T_-(g_2)) \cdot (T_+(h), br_h^\prime, T_-(h))
=(T_+(g_2), br_2 br_h^\prime, T_-(h))
\end{align*}
has branches $uz_j \to vq_j$, $j=1, \dots, n$. 
We recall that $z_j \to q_j$, $j=1, \dots, n$ 
are branches of a reduced tree-braid-tree diagram. 
Hence, 
we can not reduce the carets formed by the $u z_j \to v q_j$ 
of the tree-braid-tree diagram $(T_+(g_2), br_2 br_h^\prime, T_-(h))$. 

Since $\ell_1 \leq N(g_1)$, we have
\begin{align*}
N(g_3) = N(g_2h) 
&\geq N(x_{m-\ell_1}) \\
&=m-\ell_1+2 \\
&\geq m-N(g_1)+2 \\
&=M(N(g_1)+1)+1-N(g_1)+2 \\
&=(M-1)N(g_1) +M+3, 
\end{align*}
as required. 

For part (3), let $r=\ell_0(T_+(g_1))$. 
We note that $r \leq N(g_1)$ holds by the same reason as $\ell_1 \leq N(g_1)$. 
Let $s$ be a binary word 
such that $0^r \to s$ is a branch of $(T_+(g_1), br_1, T_-(g_1))$ of $g_1$. 
By the definition of $br_h$, 
$s \to 0^r$ is a branch of the diagram $(T_-(g_1), br_h, T_+(g_1))$ of $h$. 
Recall that $u \to 1^{\ell_1}$ is a branch of 
$(T_+(g_1), br_1, T_-(g_1))$ of $g_1$. 
We consider two cases: (a) $u \not \equiv 0^r$, and (b) $u \equiv 0^r$. 

In case (a), $0^r \to s$ is a branch of $(T_+(g_2), br_2, T_-(g_2))$. 
Indeed, every branch of $(T_+(g_1), br_1, T_-(g_1))$ of $g_1$, 
other than $u \to 1^{\ell_1}$, is also a branch of $(T_+(g_2), br_2, T_-(g_2))$. 
Then since $(T_+(g_2), br_2, T_-(g_2))$ has the branch $0^r \to s$ and 
$(T_+(g_1), br_h, T_+(g_1))$ of $h$ has the branch $s \to 0^r$, 
by adding some minimal number of caret if necessary, 
the tree-braid-tree diagram of the product $g_2h$ has the branch $0^r \to 0^r$. 
Since this diagram might be not reduced, $\ell_0(g_3) \leq r$ holds. 
Hence, by $r \leq N(g_1)$ holds, we have $\ell_0(g_3) \leq N(g_1)$ and
$0^{\ell_0(g_3)} \to 0^{\ell_0(g_3)}$ is reduced tree-braid-tree diagram of 
$g_2h=g_3$, as required. 
We illustrate a sketch of the product $g_2h$ in Figure \ref{image_of_subpath3}. 

\begin{figure}[tbp]
	\centering
	\includegraphics[width=110mm]{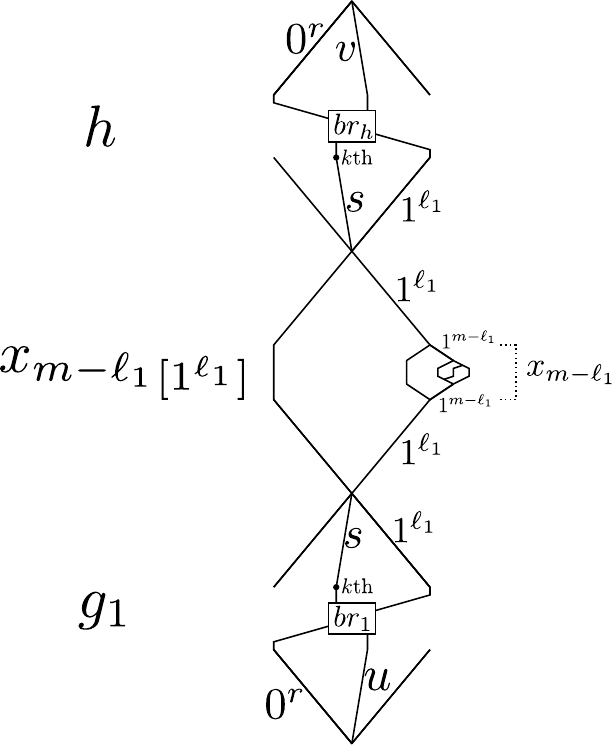}
	\caption{A rough sketch of a tree-braid-tree diagram of the product $g_1 {x_{m-\ell_1}}_{[1^{\ell_1}]}h$ with only the important branches}
	\label{image_of_subpath3}
\end{figure}
In case (b), $u \equiv 0^r$, so $1^{\ell_1} \equiv s$. 
Since the diagram $(T_-(g_1),br_h, T_+(g_1))$ of $h$ 
has the branches $1^{\ell_1} \to v$ and $s \to 0^r$, 
we have $v \equiv 0^r \equiv u$. 
Since the reduced tree-braid-tree diagram of $g_2 h=g_3$ 
has the branches $u z_j \to v q_j$, 
this diagram has the branches $0^r z_j \to 0^r q_j$. 
Since $m-\ell_1>0$ holds, 
$x_{m-\ell}$ has a branch $0 \to 0$. 
Hence, reduced tree-braid-tree diagram of $g_3$ has a branch
$0^{r+1} \to 0^{r+1}$, as required. 
\end{proof}
\end{lemma}
\begin{subpath}We fix an integer $Q \geq 12M/C_1^2$ and let
\begin{align*}
w_4  \equiv x_0^{Q|g|}x_{1}^{-1}x_0^{-Q|g|+1}. 
\end{align*}
We also let $g_4=g_3w_4$. 
We note that $w_4$ is a word representing 
the terminal vertex of Proposition $\ref{braided_main_prop}$. 
\end{subpath}
Part (1) of the following lemma corresponds to \cite[Lemma 2.11]{golan2019divergence}. 
\begin{lemma} \label{braided_subpath4}
The following assertions hold. 
\begin{enumerate}[font=\normalfont]
\item For every prefix $w^{\prime}$ of $w_4$ we have
\begin{align*}
N(g_3w^{\prime}) \geq N(g_3) 
+ \frac{1}{2} \| w^{\prime} \| -2N(g_1)-1. 
\end{align*}
\item As elements in $BV$, $g_3$ and $w_4$ commute. 
\end{enumerate}
\begin{proof}
To prove part (1), let $\ell=\ell_0(g_3)$. 
We 
first consider prefixes $w^\prime$ of $w_4$ which are positive power of $x_0$. 
We note that by Lemmas \ref{braided_subpath3} (2), 
\ref{braided_subpath1} and \ref{braided_subpath3} (3), we have
\begin{align*}
N(g_3) 
&\geq (M-1)N(g_1) + M + 3 \\
&\geq N(g_1)+N(g_1)+M+3 \\
&\geq N(g) + \ell + M + 2 \\
&\geq 3 + \ell -1 . 
\end{align*}
Hence, we can apply Corollary \ref{braided_cor_times_x0_2} (1) to $g_3$. 
Again we note that $\ell -1 \leq N(g_1)$ holds by Lemma \ref{braided_subpath3} (3). 
Then we have 
\begin{align}
N(g_3w^\prime) = N(g_3x_0^i)
&\geq N(g_3) + i -2(\ell -1) \nonumber \\
&\geq N(g_3) + i -2N(g_1) \label{braided_g3wprime_hyouka} \\
&= N(g_3) + \| w^\prime \| -2N(g_1) \nonumber \\
&\geq N(g_3) + \frac{1}{2} \| w^\prime \| -2N(g_1).  \nonumber
\end{align}
Thus, to finish the proof of part (1), it suffices to show that 
for every prefix $w^\prime$ of $w_4$ which contains the letter $x_1^{-1}$, 
we have 
\begin{align}
N(g_3w^\prime) \geq N(g_3x_0^{Q|g|}) -1. \label{braided_subpath4_goal}
\end{align}
Indeed, in that case, by inequality \eqref{braided_g3wprime_hyouka}, 
\begin{align*}
N(g_3w^\prime) 
&\geq N(g_3x_0^{Q|g|})-1 \\
&\geq N(g_3) + Q|g| -2N(g_1) -1 \\
&= N(g_3) + \frac{1}{2} \| w_4 \| - 2N(g_1)-1 \\
&\geq N(g_3) + \frac{1}{2}\| w^\prime \| -2N(g_1) -1, 
\end{align*}
as desired. 
To show inequality \eqref{braided_subpath4_goal}, 
we first consider the following prefix 
\begin{align*}
p \equiv x_0^{Q|g|}x_1^{-1}x_0^{-1} 
\equiv x_0^{Q|g|-2} \cdot x_0^2x_1^{-1}x_0^{-1}. 
\end{align*}
Since $C_1Q \geq 1200$ and $N(g) \geq 3$ hold, by Theorem \ref{braided_bound}, 
we note that we have 
\begin{align*}
\frac{1}{C_1}|g|-Q|g| \leq \frac{|g|}{C_1}-\frac{1200|g|}{C_1}
= \frac{-1199}{C_1}|g| \leq -1199N(g) < -5. 
\end{align*}
By Lemmas \ref{braided_subpath3} (3), \ref{braided_subpath1} and 
Theorem \ref{braided_bound}, we have 
\begin{align*}
\ell \leq N(g_1)+1 \leq N(g)+3 \leq \frac{1}{C_1}|g|+3 < Q|g| -2. 
\end{align*}
Since $N(g_3) \geq 3+\ell-1$ holds, by Lemma \ref{braided_lem_times_x0_2} (1) and (2)
(if necessary, apply them repeatedly), we have $\ell_0(g_3x_0^{Q|g|-2})=1$. 
In other words, 
the range-tree of the reduced tree-braid-tree diagram of $g_3x_0^{Q|g|-2}$ 
has a branch $0$. 
By applying Lemma \ref{braided_lem_times_x0_2} (1) to $g_3x_0^{Q|g|-2}$ twice, 
we have 
\begin{align*}
N(g_3x_0^{Q|g|})
=N(g_3x_0^{Q|g|-2} \cdot x_0^2)
=N(g_3x_0^{Q|g|-2} \cdot x_0)+1
=N(g_3x_0^{Q|g|-2})+2. 
\end{align*}
Since $x_0^2x_1^{-1}x_0^{-1}={x_0}_{[0]}$ and $\ell_0(g_3x_0^{Q|g|-2})=1$ hold, 
by Lemma \ref{braided_add_copy}, we have
\begin{align*}
N(g_3p)
&=N(g_3x_0^{Q|g|-2} \cdot x_0^2x_1^{-1}x_0^{-1}) \\
&=N(g_3x_0^{Q|g|-2} \cdot {x_0}_{[0]}) \\
&=N(g_3x_0^{Q|g|-2}) + N(x_0) \\
&=N(g_3x_0^{Q|g|}) -2 +2 \\
&=N(g_3x_0^{Q|g|}), 
\end{align*}
so inequality \eqref{braided_subpath4_goal} holds for the prefix $p$. 
We also note that by Lemma \ref{braided_add_copy}, $\ell_0(g_3p) \neq 1$ 
and $N(g_3p) \geq 3$. 

To finish the proof, it remains to prove that inequality \eqref{braided_subpath4_goal} 
holds 
for the prefix (a) $w^\prime \equiv x_0^{Q|g|}x_1^{-1}$ 
and prefixes of the form 
(b) $w^\prime \equiv x_0^{Q|g|}x_1^{-1}x_0^{-i}$, for $1 < i \leq Q|g|-1$. 
For the case (a), we note that $\ell_0(g_3p) \neq 1$ 
and $g_3w^\prime = g_3px_0$. 
Hence, by applying Lemma \ref{braided_lem_times_x0_2} (2) 
to $g_3p$, we have 
\begin{align*}
N(g_3w^\prime) = N(g_3px_0) \geq N(g_3p)-1=N(g_3x_0^{Q|g|})-1, 
\end{align*}
as required. 
Finally, we note again that $\ell_0(g_3p) \neq 1$ and 
prefixes of the form (b) can be written as  $w^\prime \equiv px_0^{-(i-1)}$. 
Hence, by Corollary \ref{braided_cor_times_x0^-1_2} (2) we have
\begin{align*}
N(g_3w^\prime)=N(g_3px_0^{-(i-1)}) \geq N(g_3p) = N(g_3x_0^{Q|g|}), 
\end{align*}
as required. 

\begin{figure}[tbp]
	\centering
	\includegraphics[width=80mm]{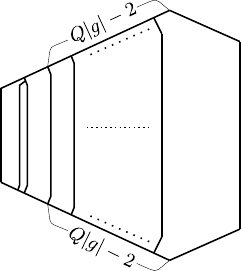}
	\caption{A tree-braid-tree diagram of $w_4$}
	\label{w4_diagram}
\end{figure}
For part (2), we first note that reduced tree-braid-tree diagram of $g_3$ 
has a ``same length'' branch $0^\ell \to 0^\ell$. 
By calculating $x_0^i (x_0^2x_1^{-1}x_0^{-1})x_0^{-i}$, 
for $i=1, 2, \dots, Q|g|-2$ inductively, 
we get tree-braid-tree diagram of $w_4$ as Figure \ref{w4_diagram}. 
Since $\ell-1 <Q|g|-2$ holds, we can calculate $w_4g_3$ and $g_3w_4$ 
as Figure \ref{w4_and_g3_kakan} and Figure \ref{g3_and_w4_kakan}, respectively. 
Hence, we have $w_4g_3=g_3w_4$, as required. 
\begin{figure}[tbp]
	\centering
	\includegraphics[width=110mm]{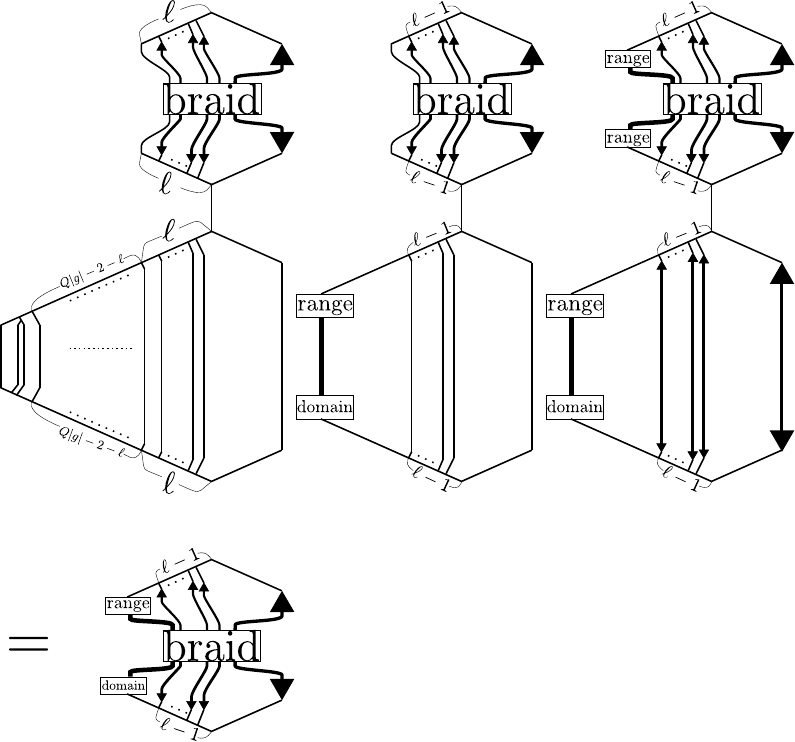}
	\caption{Calculation of $w_4g_3$}
	\label{w4_and_g3_kakan}
	\centering
	\includegraphics[width=110mm]{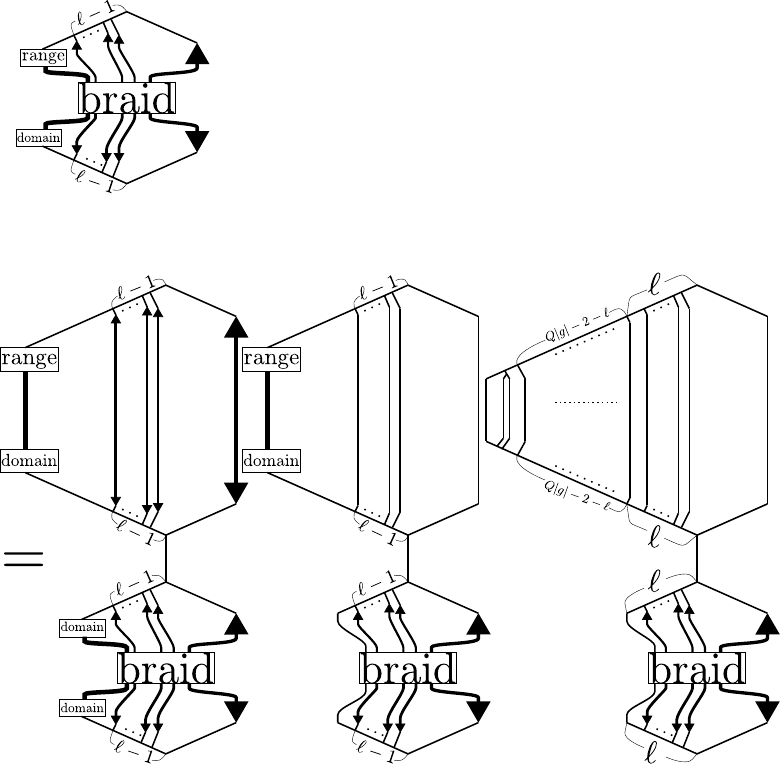}
	\caption{Calculation of $g_3w_4$}
	\label{g3_and_w4_kakan}
\end{figure}
\end{proof}
\end{lemma}
\begin{subpath}
Let $w_5$ be a minimal word in the alphabet $X$ 
such that $w_5=g_3^{-1}$ in $BV$. 
Let $g_5=g_4w_5$. 
\end{subpath}
It follows from Lemma \ref{braided_subpath4} (2) that 
\begin{align*}
gw=gw_1w_2w_3w_4w_5=g_5=g_3w_4g_3^{-1}=w_4. 
\end{align*}
Hence, $gw=x_0^{Q|g|}x_1^{-1}x_0^{-Q|g|+1}$ for $w \equiv w_1w_2w_3w_4w_5$, 
as required. 

It remains to prove that one can choose constants $\delta$, $D$ (independently of $g$), 
so that path $w$ satisfies the conditions in the Proposition \ref{braided_main_prop}. 
First, by definitions of subpaths, we have the following. 
\begin{align*}
\| w \|
&\leq \| w_1w_2w_3 \| + \| w_4 \| + \| w_5 \| \\
&= \| w_1w_2w_3 \| + 2Q|g| + |g_3| \\
&= \| w_1w_2w_3 \| + 2Q|g| + |gw_1w_2w_3| \\
&\leq \| w_1w_2w_3 \| + 2Q|g| + |g| + \| w_1w_2w_3 \| \\
&= 2\| w_1w_2w_3 \| + (2Q+1)|g| \\
&\leq 2\| w_1w_2w_3 \| + 3Q|g|. 
\end{align*}
Furthermore, we have a upper bound of $\| w_1w_2w_3 \|$ as follows. 
\begin{align}
\| w_1w_2w_3 \|
&\leq \| w_1 \| + \| w_2 \| + \| w_3 \| \nonumber \\
&\leq 4 + 2M(N(g_1)+1)+1 + 14N(g_1) \nonumber \\
&=2MN(g_1) + 14N(g_1) + 5 + 2M \nonumber \\
&\leq 2M(N(g)+2) + 14(N(g)+2) + 5 + 2M \nonumber \\
&= 2MN(g) + 14N(g) + 33 + 6M \nonumber \\
&< 2MN(g) + 14N(g) + 33N(g) + 2M \times 3 \nonumber \\
&=2MN(g) + 47N(g) + 2M \times 3 \nonumber \\
&\leq 2MN(g) + MN(g) + 2MN(g) \nonumber \\
&= 5MN(g) \label{braided_w1w2w3_and_N(g)} \\
&\leq \frac{5M}{C_1}|g|,  \tag*{$(4.7)'$} \label{braided_w1w2w3_and_|g|)}
\end{align}
where these inequalities follow from  the definition of the subpaths, 
Lemmas \ref{braided_subpath3} (1), \ref{braided_subpath1}, the definition of $M$ 
and Theorem \ref{braided_bound}. 
Therefore, we have $\| w \| \leq D|g|$ where $D=10M/C_1+3Q$, as required. 
Now, let $\delta=C_1/10M$. 
The following lemma corresponds to \cite[Lemma 2.12]{golan2019divergence} and 
completes the proof of Proposition \ref{braided_main_prop}. 
\begin{lemma}
Let $w^\prime$ be a prefix of w. 
Then $|gw^\prime| > \delta|g|$. 
\begin{proof}
First, we note that by Lemma \ref{braided_subpath3} (2), 
\begin{align*}
N(g_3) \geq (M-1)N(g_1) + M + 3. 
\end{align*}
Then for each prefix $\tilde{w} \leq w_4$ we have by Lemma \ref{braided_subpath4}, 
\begin{align}
N(g_3\tilde{w}) 
&\geq N(g_3) + \frac{1}{2}\| \tilde{w} \| -2N(g_1) -1 \nonumber \\
&\geq \frac{1}{2}\| \tilde{w} \| + (M-3)N(g_1)+M+2  \nonumber \\
&> \frac{1}{2}\| \tilde{w} \| \label{braided_g3wtilde_wtilde}. 
\end{align}

We separate the proof into two cases depending on the length of $g$. 

Case (1): $|g| < 10MN(g)$. 

It follows from Lemmas \ref{braided_subpath1} and \ref{braided_subpath2} (1) that
for every prefix $w^\prime \leq w_1w_2$, we have $N(gw^\prime) \geq N(g)$. 
Then, by applying Theorem \ref{braided_bound} to $gw^\prime$, we have 
\begin{align*}
|gw^\prime| \geq C_1N(gw^\prime)
\geq C_1N(g)
> \frac{C_1}{10M}|g|
=\delta |g|, 
\end{align*}
as required. 
Next, we consider a prefix $w^\prime \leq w_3$. 
By Theorem \ref{braided_bound}, Lemma \ref{braided_subpath2} (2) and 
$M \geq 100/C_1$, 
\begin{align*}
|g_2| \geq C_1N(g_2) \geq C_1MN(g_1) \geq 100N(g_1). 
\end{align*}
Since we already know that 
$\| w_3 \| \leq 14N(g)$ (Lemma \ref{braided_subpath3} (1)), 
$N(g_1) \geq N(g)$ (Lemma \ref{braided_subpath1}) and
$N(g) > |g|/10M$ (assumption of case (1)) hold, we have
\begin{align*}
|g_2w^\prime|
&\geq |g_2|-\| w^\prime \| \\
&\geq |g_2| - \| w_3 \| \\
&\geq 100N(g_1) - 14N(g_1) \\
&= 86N(g_1) \\
&\geq 86N(g) \\
&> \frac{86}{10M}|g| \\
&> \frac{C_1}{10M}|g|  = \delta|g|, 
\end{align*}
as required.  
Now, let $w^\prime$ be a prefix of $w_4$. By Theorem \ref{braided_bound}, Lemmas \ref{braided_subpath4} (1), \ref{braided_subpath3} (2), \ref{braided_subpath1}, and assumption of case (1), we have
\begin{align*}
|g_3w^\prime| &\geq C_1N(g_3w^\prime) \\
&\geq C_1(N(g_3)-2N(g_1)-1) \\
&\geq C_1((M-1)N(g_1)+M+3-2N(g_1)-1) \\
&= C_1((M-3)N(g_1)+M+2) \\
&>C_1N(g) \\
&>\frac{C_1}{10M}|g|=\delta|g|, 
\end{align*}
as required. 
Finally, we consider a prefix $w^\prime \leq w_5$. 
Since $\| w_1w_2w_3 \| \leq (5M/C_1)|g|$ 
(inequality \ref{braided_w1w2w3_and_|g|)}) 
and $Q \geq 12M/C_1^2$, we have
\begin{align*}
\| w_5 \| = |g_3| = |gw_1w_2w_3| 
&\leq |g| + \| w_1w_2w_3 \| \\
&\leq |g| + \frac{5M}{C_1}|g| \\
&< \frac{M}{C_1}|g| + \frac{5M}{C_1}|g| \\
&= \frac{6M}{C_1}|g| \\
&\leq \frac{C_1Q}{2}|g|. 
\end{align*}
By Theorem \ref{braided_bound}, 
inequality \eqref{braided_g3wtilde_wtilde} 
for $\tilde{w}=w_4$ (then $\| w_4 \|=2Q|g|$) and 
the definition of $Q$, we have 
\begin{align*}
|g_4w^\prime| \geq |g_4| - \| w^\prime\| 
&\geq |g_4| - \| w_5 \| \\
&> C_1N(g_4) - \frac{C_1Q}{2}|g| \\
&\geq C_1Q|g| - \frac{C_1Q}{2}|g| \\
&= \frac{C_1Q}{2}|g| \\
&> \frac{C_1}{10M}|g|=\delta|g|, 
\end{align*}
as required. 
Hence, the lemma holds in case (1). 

Case (2): $|g| \geq 10MN(g)$. 

Since $\| w_1w_2w_3 \| \leq 5MN(g)$ (inequality \eqref{braided_w1w2w3_and_N(g)}) 
and $|g|/2 \geq 5MN(g)$ (assumption of case (2)), 
for any prefix $w^\prime \leq w_1w_2w_3$ we have 
\begin{align*}
|gw^\prime| 
&\geq |g| - \| w^\prime \| \\
&\geq |g| - \| w_1w_2w_3 \| \\
&\geq |g| - 5MN(g) \\
&\geq |g| - \frac{1}{2}|g| \\
&=\frac{1}{2}|g| \\
&> \frac{C_1}{10M}|g| =\delta|g|, 
\end{align*}
as required. 
In particular, we note that $|g_3| \geq |g|/2$ holds where $g_3=gw_1w_2w_3$. 
Let  $w^\prime \leq w_4$. 
If $\| w^\prime \| \leq |g|/5$, then we have
\begin{align*}
|g_3w^\prime| 
\geq |g_3|-\| w^\prime \|
\geq \frac{1}{2}|g| - \frac{1}{5}|g|
= \frac{3}{10}|g|
> \frac{C_1}{10M}|g|
=\delta |g|, 
\end{align*}
as desired. 
Hence, we can assume that $\| w^\prime \| > |g|/5$. 
In that case, by Theorem \ref{braided_bound} and 
inequality \eqref{braided_g3wtilde_wtilde} we have 
\begin{align*}
|g_3w^\prime| 
\geq C_1N(g_3w^\prime)
> \frac{C_1}{2}\| w^\prime \|
> \frac{C_1}{10}|g|
> \frac{C_1}{10M}|g|
=\delta |g|, 
\end{align*}
as required. 
To finish the proof, we note that by Theorem \ref{braided_bound} and 
inequality \eqref{braided_g3wtilde_wtilde} 
for $\tilde{w}=w_4$ (then $\| w_4 \|=2Q|g|$), we have 
\begin{align*}
|g_4| = |g_3w_4| 
\geq C_1N(g_3w_4)
> \frac{C_1}{2} \| w_4 \|
=C_1Q|g|. 
\end{align*}
By inequality \eqref{braided_w1w2w3_and_N(g)} and assumption of case (2), 
we also note that we have 
\begin{align*}
\| w_5 \| = |g_3| 
&= |gw_1w_2w_3| \\
&\leq |g| + \| w_1w_2w_3 \| \\
&\leq |g| + 5MN(g) \\
&\leq |g| + \frac{1}{2}|g| \\
&=\frac{3}{2}|g|. 
\end{align*}
Hence, since $Q \geq 12M/C_1^2$, we have
\begin{align*}
|g_4w^\prime| \geq |g_4| - \| w^\prime \| 
&\geq |g_4| - \| w_5 \| \\
&> C_1Q|g| - \frac{3}{2}|g| \\
&= (C_1Q - \frac{3}{2})|g| \\
&> |g| \\
&> \frac{C_1}{10M}|g| = \delta |g|, 
\end{align*}
as required. 
\end{proof}
\end{lemma}
\end{proof} 
\end{proposition}

\begin{figure}[tbp]
	\centering
	\includegraphics[width=150mm]{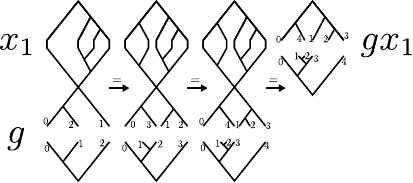}
	\caption{Calculating example of $gx_1$}
	\label{1-1}
\end{figure}

\begin{figure}[tbp]
	\centering
	\includegraphics[width=150mm]{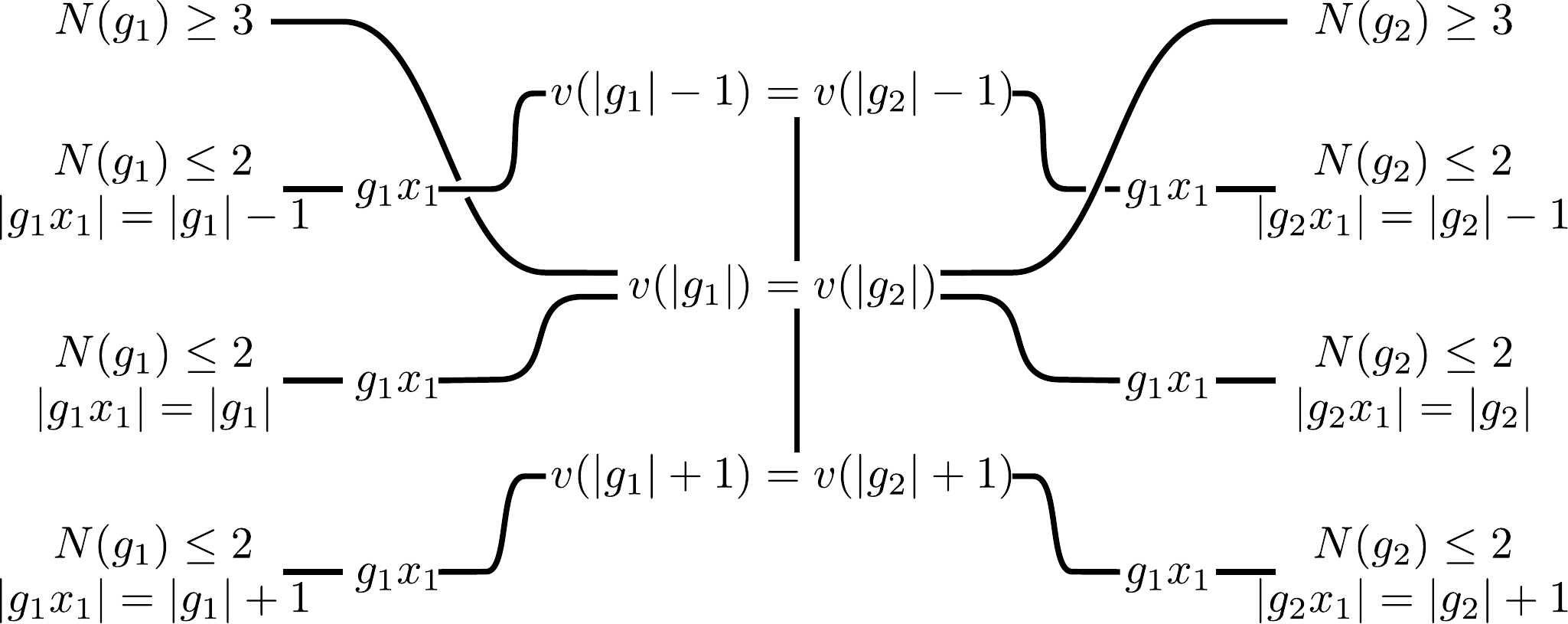}
	\caption{The path connecting $g_1$ and $g_2$ where $|g_1|, |g_2| \geq 2$}
	\label{chart}
\end{figure}

By multiplying $x_1$ on the right to an element having one or two carets, 
we get  the element of $BV$ 
satisfying the assumption of Proposition \ref{braided_main_prop}. 
\begin{lemma}  \label{times_x1_BV}
Let $g \in BV$ be such that $N(g) \leq 2$. Then $N(gx_1) \geq 3$. 
\begin{proof}
By regarding each braid as just a permutation, 
it can be shown by finite number of direct calculations. 
Indeed, each tree-braid-tree diagram is reduced if 
there exists no strands pair such that they have a same parent. 
Hence, if it is reduced when considering $g$ as the element of $V$, 
then it is also reduced in $BV$. 
For example, see Figure \ref{1-1}. 
The endpoints of each strand are represented by the same number, 
with a blank representing some braid. 
\end{proof}
\end{lemma}
\begin{lemma}\label{bounded_-1}
Let $g \in BV$. 
Then
\begin{align*}
|g|-1 \leq |gx_1| \leq |g|+1. 
\end{align*}
\begin{proof}
The first inequality follows from $|g| \leq |gx_1|+|x_1^{-1}|=|gx_1|+1$. 
The second inequality follows from $|gx_1| \leq |g| + |x_1| = |g|+1$. 
\end{proof}
\end{lemma}
The following proposition immediately implies that 
braided Thompson group has liner divergence, 
completing the proof of Theorem \ref{main_theorem_BV}. 
See Figure \ref{chart} for the overview of the paths. 
The idea of paths (corresponding two vertical lines in the middle) in the following proposition 
comes from \cite[Theorem 2.13]{golan2019divergence}. 
\begin{proposition} \label{braided_sub_prop}
There exist constants $\delta_{BV}$ and $D_{BV} > 0$ such that the following holds. 
Let $g \in BV$ be an element with $|g| \geq 2$. 
Then there exists a path of length at most $D_{BV}|g|$ in the Cayley graph 
$\Gamma = \mathrm{Cay}(BV, X)$ which avoids a $\delta_{BV}|g|$-neighborhood 
of the identity and which has initial vertex $g$ 
and terminal vertex $x_0^{Q|g|}x_1^{-1}x_0^{-Q|g|+1}$. 

In other words, there exists a word $w_{BV}$ in the alphabet $X$ 
such that $ \|w \| < D_{BV}|g|$; for any prefix $w^{\prime}$ of $w_{BV}$, 
we have $|gw^{\prime}| > \delta_{BV}|g|$ and such that 
\begin{align*}
gw_{BV} = x_0^{Q|g|}x_1^{-1}x_0^{-Q|g|+1}. 
\end{align*}
\begin{proof}
For each natural number $k >0$, let 
\begin{align*}
v(k)=x_0^{Q k}x_1^{-1}x_0^{-Q k+1}. 
\end{align*}
First, if $N(g) \geq 3$ then the proposition follows from 
Proposition \ref{braided_main_prop}. 
Hence, we can assume that $N(g) \leq 2$. 
Then by Lemma \ref{bounded_-1}, $|gx_1|=|g|-1$, $|g|$ or $|g|+1$ and
by Lemma \ref{times_x1_BV}, $N(gx_1) \geq 3$. 
Let 
\begin{align*}
D_{BV } &=2D+4Q+1 \\
\delta_{BV} &= \min \{ \frac{1}{2}\delta, \frac{1}{2}C_1Q \}. 
\end{align*}
In the following, we will use Proposition \ref{braided_main_prop} 
to construct the path connecting $gx_1$ and $v(|g|)$. 
We consider three cases depending on the length of $gx_1$. 

Case (1): $|gx_1|=|g|-1$. 

By Proposition \ref{braided_main_prop}, 
there exists a path of length at most $D(|g|-1)$ 
which avoids a $\delta(|g|-1)$-neighborhood of identity and 
which has initial vertex $gx_1$ and terminal  vertex $v(|g|-1)$. 
Since $|g| \geq 2$, we have $\delta(|g|-1) \geq (\delta /2)|g|$. 
Hence this path avoids a $(\delta/2)|g|$-neighborhood of identity. 
Thus, we construct a path connecting $v(|g|-1)$ and $v(|g|)$. 
Let 
\begin{align*}
p(|g|-1) \equiv x_0^{Q(|g|-1)-1}x_1x_0^Q x_1^{-1}x_0^{-Q|g|+1}. 
\end{align*}
It is clear that $p(|g|-1)$ labels a path from $v(|g|-1)$ to $v(|g|)$ and 
the length of $p(|g|-1)$ is at most $2Q|g|$. 
In the following, we prove that for any prefix $p^\prime$ of $p(|g|-1)$, 
we have $|v(|g|-1)p^\prime| > C_1Q|g|$. 
Indeed, it is easy to see that the positive part of the normal form of 
the element $v(|g|-1)p^\prime$ is $x_0^i$ for $i \geq Q(|g|-1)$. 
Hence, by \cite[Theorem 3]{burillo2001metrics}, we have
\begin{align*}
N(v(|g|-1)p^\prime)
\geq N(x_0^i)
\geq Q(|g|-1)+1
> Q(|g|-1), 
\end{align*}
where we note that this theorem claims only the relationship between 
the number of carets of elements in $F$ and their exponents, 
so it can be applied to $BV$. 
Hence, by Theorem \ref{braided_bound}, we have
\begin{align*}
|v(|g|-1)p^\prime| \geq C_1N(v(|g|-1)p^\prime) > C_1Q(|g|-1). 
\end{align*}
Since $|g| \geq 2$, we have $C_1Q(|g|-1) \geq (C_1Q/2)|g|$, as required. 

Case (2): $|gx_1|=|g|$. 

By Proposition \ref{braided_main_prop}, all assertions follow. 

Case (3): $|gx_1|=|g|+1$. 

By Proposition \ref{braided_main_prop}, 
there exists a path of length at most $D(|g|+1)$ 
which avoids a $\delta(|g|+1)$-neighborhood of identity and 
which has initial vertex $gx_1$ and terminal  vertex $v(|g|+1)$. 
Since $|g| \geq 2$, we have $D(|g|+1) \leq 2D|g|$. 
Thus, we construct a path connecting $v(|g|)$ and $v(|g|+1)$. 
Let 
\begin{align*}
p(|g|) \equiv x_0^{Q(|g|)-1}x_1x_0^Q x_1^{-1}x_0^{-Q(|g|+1)+1}. 
\end{align*}
It is clear that $p(|g|)$ labels a path from $v(|g|)$ to $v(|g|+1)$ and 
the length of $p(|g|)$ is at most $2Q(|g|+1)$. 
Since $|g| \geq 2$, we have $2Q(|g|+1) \leq 4Q|g|$. 
By the almost same argument as case (1), we have 
\begin{align*}
|v(|g|)p^\prime| \geq C_1N(v(|g|)p^\prime) 
\geq C_1N(x_0^i) > C_1Q|g|, 
\end{align*}
for any prefix $p^\prime \leq p(|g|)$ and corresponding $i \geq Q|g|$, 
as required. 
\end{proof}
\end{proposition}
\bibliographystyle{plain}
\bibliography{hoge} 
\address{
Department of Mathematical Sciences, 
Tokyo Metropolitan University,
Minami-osawa Hachioji, Tokyo, 192-0397, Japan
}

\textit{E-mail address}: \href{mailto:kodama-yuya@ed.tmu.ac.jp}{\texttt{kodama-yuya@ed.tmu.ac.jp}}
\end{document}